\documentclass[a4paper,12pt,english]{amsart}

\usepackage{amsmath,amssymb,amsthm,graphicx,pstricks,xypic,comment,fancyhdr}
\usepackage[all]{xy}
\usepackage{epic}
\usepackage{eufrak}
\usepackage{lscape}
\usepackage[english]{babel}

\title[Cluster categories: a generalization]{Cluster categories for algebras of global dimension 2 and quivers with potential}
\author{Claire Amiot}
\address{Universit{\'e} Paris 7\\
Institut de Math{\'e}matiques de Jussieu \\
Th{\'e}orie des groupes et des repr{\'e}sentations\\
Case 7012 \\ 2 Place Jussieu \\
75251 Paris Cedex 05, France}
\email{amiot@math.jussieu.fr}

\newcommand{\Hom}{{\sf Hom }}
\newcommand{\End}{{\sf End }}
\newcommand{\Ext}{{\sf Ext }}
\newcommand{\Tor}{{\sf Tor}}
\newcommand{\Ker}{{\sf Ker }}
\newcommand{\Coker}{{\sf Coker }}
\newcommand{\Imm}{{\sf Im }}
\renewcommand{\mod}{{\sf mod \hspace{.02in}  }}
\newcommand{\Mod}{{\sf Mod \hspace{.02in} }}

\newcommand{\ind}{{\sf ind \hspace{.02in} }}
\newcommand{\per}{{\sf per \hspace{.02in}  }}
\newcommand{\proj}{{\sf proj \hspace{.02in} }}

\newcommand{\Sub}{{\sf Sub \hspace{.02in} }}

\newcommand{\ten}{\otimes}
\newcommand{\lten}{\overset{\boldmath{L}}{\ten}}
\newcommand{\limproj}{\underset{\leftarrow}{\rm{lim}}}

\newcommand{\Bb}{\mathcal{B}}
\newcommand{\Cc}{\mathcal{C}}
\newcommand{\Dd}{\mathcal{D}}

\newcommand{\Ff}{\mathcal{F}}
\newcommand{\Hh}{\mathcal{H}}
\newcommand{\Ii}{\mathcal{I}}
\newcommand{\Jj}{\mathcal{J}}
\newcommand{\Tt}{\mathcal{T}}
\newcommand{\Mm}{\mathcal{M}}
\newcommand{\Nn}{\mathcal{N}}
\newcommand{\Pp}{\mathcal{P}}

\newcommand{\Mmm}{\overline{\mathcal{M}}}

\newcommand{\bsm}{\begin{smallmatrix}}
\newcommand{\esm}{\end{smallmatrix}}

\newtheorem{thma}{Theorem}[section]
\newtheorem*{thm*}{Th{\'e}or{\`e}me}
\newtheorem{lema}[thma]{Lemma}
\newtheorem*{lem*}{Lemme}

\newtheorem{cora}[thma]{Corollary}
\newtheorem*{prop*}{Proposition}
\newtheorem{prop}[thma]{Proposition}
\theoremstyle{remark}

\theoremstyle{definition}

\newtheorem{dfa}[thma]{Definition}

\begin{document}

\maketitle

\begin{abstract}
Let $k$ be a field and $A$ a finite-dimensional $k$-algebra of global dimension $\leq 2$.
We construct a triangulated category $\Cc_A$ associated to $A$ which, if $A$ is hereditary, is triangle equivalent to the cluster category of $A$. When $\Cc_A$ is $\Hom$-finite, we prove that it is 2-CY and endowed with a canonical cluster-tilting object. This new class of categories contains some of the stable categories of modules over a preprojective algebra studied by Geiss-Leclerc-Schr{\"o}er and by Buan-Iyama-Reiten-Scott. Our results also apply to quivers with potential. Namely, we introduce a cluster category $\Cc_{(Q,W)}$ associated to a quiver with potential $(Q,W)$. When it is Jacobi-finite we prove that it is endowed with a cluster-tilting object whose endomorphism algebra is isomorphic to the Jacobian algebra $\Jj(Q,W)$. 
 
\end{abstract}
\tableofcontents

\section*{Introduction}

The cluster category associated with a finite-dimensional hereditary algebra was introduced in \cite{Bua} (and in \cite{Cal} for the $A_n$ case). It serves in the representation-theoretic approach to cluster algebras introduced and  studied by Fomin and Zelevinsky in a series of articles (cf. \cite{FZ1}, \cite{FZ2}, \cite{FZ4} and \cite{BFZ} with Berenstein). The link between cluster algebras and cluster categories is in the spirit of `categorification'. Several articles (e.g. \cite{Mar}, \cite{Bua},
 \cite{Cal2}, \cite{CC}, \cite{BMR1}, \cite{BMR2},
\cite{BMRT}, \cite{CK2}) deal with the categorification of the cluster algebra $A_Q$ associated with an acyclic quiver $Q$ using the cluster category $\Cc_Q$ associated with the path algebra of the quiver $Q$.
Another approach consists in categorifying cluster algebras by subcategories of the category of modules over a preprojective algebra associated to an acyclic quiver (cf. \cite{Gei}, \cite{Gei2}, \cite{Gei3},
\cite{Gei4}, \cite{Bua2}). In both approaches the categories $\Cc$ (or their associated stable categories) satisfy the following fundamental properties:

\hspace{.2in}- $\Cc$ is a triangulated category;

\hspace{.2in}- $\Cc$ is $2$-Calabi-Yau (2-CY for short); 
  
\hspace{.2in}- there exist cluster-tilting objects.\\
It has been shown that these properties alone imply many of the most important theorems about cluster categories, respectively stable module categories over preprojective algebras
  (cf. \cite{Iya}, \cite{Kel4},
 \cite{Kel5}, \cite{Kel7}, \cite{Pal}, \cite{Tab}). In particular by \cite{Iya}, in a category $\Cc$ with such properties it is possible to `mutate' the cluster-tilting objects and there exist exchange triangles. This is fundamental for categorification.

Let $k$ be a field. In this article we want to generalize the construction of the cluster category replacing the hereditary algebra $kQ$ by a finite-dimensional algebra $A$ of finite global dimension.  A candidate might be the orbit category
$\Dd^b(A)/\nu [-2],$ where $\nu$ is the Serre functor of the derived
category $\Dd^b(A)$. By \cite{Kel}, such a category is
triangulated if $A$ is derived equivalent to an hereditary category
$\Hh$. However in general, it is not triangulated. 
Thus a more appropiate candidate is the triangulated hull $\Cc_A$ of
the orbit category $\Dd^b(A)/\nu [-2]$. It is defined in \cite{Kel} as
 the stabilization of a certain dg category and contains the orbit category as a full subcategory. More precisely the category $\Cc_A$  is a quotient of a triangulated category $\Tt$ by a thick subcategory $\Nn$ which is $3$-CY.  This leads us to the study of such quotients in full generality. We prove that the quotient is $2$-CY if the objects of $\Tt$ are `limits' of objects of $\Nn$ (Theorem~\ref{nondegenere}). In particular we deduce that 
the cluster category $\Cc_A$ of an algebra of finite global dimension is $2$-CY if it is $\Hom$-finite (Corollary~\ref{2CY}). 

We study the particular case where the algebra is of global dimension $\leq 2$.  Since $kQ$ is  a cluster-tilting object of the category $\Cc_Q$, the canonical candidate to be a cluster-tilting object in the category $\Cc_A$ would be $A$ itself. Using generalized tilting theory (cf. \cite{Kel9}), we
give another construction of the cluster category. We find a triangle equivalence $$\xymatrix{\Cc_A\ar[r]^-\sim & \per \Pi/\Dd^b\Pi}$$ where $\Pi$ is a dg algebra in negative degrees which is bimodule $3$-CY and homologically smooth. This equivalence sends the object $A$ onto the image of the free dg module $\Pi$ in the quotient. This leads us to the study of the categories $\per \Gamma/\Dd^b\Gamma$ where $\Gamma$ is a dg algebra with the above properties. We prove that if the zeroth cohomology of $\Gamma$ is finite-dimensional, then the category $\per \Gamma/\Dd^b\Gamma$ is $2$-CY and the image of the free dg module $\Gamma$ is a cluster-tilting object (Theorem \ref{clustertilting}). We show that the algebra $H^0\Gamma$ is finite-dimensional if and only if the quotient $\per \Gamma/\Dd^b\Gamma$ is $\Hom$-finite.  Thus we prove the existence of a cluster-tilting object in cluster categories $\Cc_A$ associated with algebras of global dimension $2$ which are $\Hom$-finite (Theorem~\ref{Aclustertilting}). Moreover, this general approach applies to the Ginzburg dg algebras \cite{Gin} associated with a quiver with potential. Therefore we introduce a new class of $2$-CY categories $\Cc_{(Q,W)}$ endowed with a cluster-tilting object associated with a Jacobi-finite quiver with potential $(Q,W)$ (Theorem~\ref{potentiel2CY}). 

In \cite{Gei4}, Geiss, Leclerc and Schr{\"o}er construct 
subcategories $\Cc_M$ of $\mod \Lambda$ (where $\Lambda=\Lambda_Q$ is a
preprojective algebra of an acyclic quiver) associated with
certain terminal $kQ$-modules $M$. We show in the last part that the stable
category of such a Frobenius category $\Cc_M$ is triangle equivalent to
a cluster category $\Cc_A$ where $A$ is the endomorphism algebra of a
postprojective module over an hereditary algebra (Theorem~\ref{casGLS}).
Another approach is given by Buan, Iyama, Reiten and Scott
in \cite{Bua2}. They construct $2$-Calabi-Yau triangulated categories
$\underline{\Sub}\Lambda/\Ii_w$ where $\Ii_w$ is a two-sided ideal of
the preprojective algebra $\Lambda=\Lambda_Q$ associated with an element
$w$ of the Weyl group of $Q$. For certain elements $w$ of the Weyl
group (namely those coming from preinjective tilting modules), we
construct a triangle equivalence between
$\underline{\Sub}\Lambda/\Ii_w$ and a cluster category $\Cc_A$ where
$A$ is the endomorphism algebra of a postprojective module over a
concealed algebra (Theorem~\ref{casBIRS}).

\subsection*{Plan of the paper}

The first section of this paper is devoted to the study of Serre functors in quotients of triangulated categories. In Section 2, we prove the existence of a cluster-tilting object in a $2$-CY category coming from a bimodule $3$-CY dg algebra. Section 3 is a direct application of these results to Ginzburg dg algebras associated with quivers with potential. In particular we give the definition of the cluster category $\Cc_{(Q,W)}$ of a Jacobi-finite quiver with potential $(Q,W)$. In section 4 we define cluster categories of algebras of finite global dimension. We apply the results of Sections 1 and 2 in subsection 4.3 to the particular case of global dimension $\leq 2$. The last section links the categories introduced in \cite{Gei4} and in \cite{Bua2} with these new cluster categories $\Cc_A$.

\subsection*{Acknowledgements}
This article is part of my Ph. D. thesis under the supervision of
Bernhard Keller. I deeply thank him for his patience and
availability. I thank Bernard Leclerc, Yann Palu
and Jan Schr\"oer for interesting and helpful discussions and Idun
Reiten for kindly answering my questions. I also would like to thank
the referee for his interesting comments and remarks.  

\subsection*{Notations}
Throughout let $k$ be a field. By triangulated category we mean $k$-linear triangulated category satisfying the Krull-Schmidt property. For all triangulated categories, we will denote the shift functor by $[1]$. For a finite-dimensional $k$-algebra $A$ we denote by $\mod A$ the category of finite-dimensional right $A$-modules. More generally, for an additive $k$-category $\Mm$ we denote by $\mod \Mm$ the category of finitely presented functors $\Mm^{op}\rightarrow \mod k$. Let $D$ be the usual duality $\Hom_k(?,k)$.
If $A$ is a differential graded (=dg) $k$-algebra, we will denote by $\Dd=\Dd A$ the derived category of dg $A$-modules and by $\Dd^b A$ its full subcategory formed by the dg $A$-modules whose homology is of finite total dimension over $k$. We write $\per A$ for the category of perfect dg $A$-modules, \emph{i.e.} the smallest triangulated subcategory of $\Dd A$ stable under taking direct summands and which contains $A$. 

\section{Construction of a Serre functor in a quotient category}

\subsection{Bilinear form in a quotient category}

Let $\Tt$ be a triangulated category and $\Nn$ a thick subcategory of $\Tt$ (\emph{i.e.} a triangulated subcategory stable under taking direct summands).  We assume that there is an auto-equivalence $\nu$ in $\Tt$ such that $\nu(\Nn)\subset \Nn$. Moreover we assume that there is a non degenerate bilinear form:

$$\beta_{N,X}: \Tt(N,X)\times \Tt(X,\nu N)\longrightarrow k$$
which is bifunctorial in $N\in \Nn$ and $X\in \Tt$.

\subsubsection*{Construction of a bilinear form in $\Tt/\Nn$}

Let $X$ and $Y$ be objects in $\Tt$. The aim of this section is to construct a bifunctorial bilinear form:
$$\beta'_{X,Y}: \Tt/\Nn(X,Y)\times \Tt/\Nn(Y,\nu X[-1]) \longrightarrow k.$$
We use the calculus of left fractions \cite{Ver} in the triangle quotient $\Tt/\Nn$. Let $s^{-1}\circ f: X\rightarrow Y$ and $t^{-1}\circ g: Y\rightarrow \nu X[-1]$ be two morphisms in $\Tt/\Nn$. We can construct a diagram
$$\xymatrix@-1.2pc{X\ar[dr]_f & & Y\ar[dl]^s\ar[dr]_g & & \nu X[-1]\ar[dl]^t \ar@(dr,r)[ddll]^{ \nu u[-1]}\\ & Y'\ar[dr] & & \nu X' [-1]\ar[dl]_{s'} & \\ & & \nu X'' [-1] & &}$$
where the cone of $s'$ is isomorphic to the cone of $s$. 
Denote by $N[1]$ the cone of $u$. It is in $\Nn$ since $\Nn$ is $\nu$-stable. Thus we get a diagram of the form:
$$\xymatrix{N\ar[r] \ar@(d,l)[dr]^v & X\ar[r]^u\ar[d]^f & X''\ar[r] & N[1] \\ & Y'\ar[d]\ar@(r,u)[dr]^w & & \\ \nu X[-1]\ar[r]_{\nu u [-1]} 
& \nu X'' [-1]\ar[r] & \nu N\ar[r] & \nu X,}$$
where the two horizontal rows are triangles of $\Tt$. We define then  $\beta'_{X,Y}$ as follows:
$$\beta'_{X,Y}(s^{-1}\circ f, t^{-1}\circ g)= \beta_{N,Y'}(v,w).$$

\begin{lema}
The form $\beta'$ is well-defined, bilinear and bifunctorial.
\end{lema}

\begin{proof}
It is not hard to check that $\beta'$ is well-defined (cf. \cite{Ami2}).
Since $\beta$ is bifunctorial and bilinear,  $\beta'$ satisfies the same properties.  
\end{proof}

\subsection{Non-degeneracy}

In this section, we find conditions on $X$ and $Y$ such that the bilinear form $\beta'_{XY}$ is non-degenerate.

\begin{dfa}
Let $X$ and $Y$ be objects in $\Tt$. A morphism $p:N\rightarrow X$ is called a \emph{local $\Nn$-cover of $X$ relative to $Y$} if $N$ is in $\Nn$ and if it induces an exact sequence:
$$\xymatrix{0\ar[r] & \Tt(X,Y)\ar[r]^{p^*} & \Tt(N,Y). }$$ 

Let $Y$ and $Z$ be objects in $\Tt$. A morphism $i:Z\rightarrow N'$ is called a \emph{local $\Nn$-envelope of $Z$ relative to $Y$} if $N'$ is in $\Nn$ and if it induces an exact sequence:
$$\xymatrix{0\ar[r] & \Tt(Y,Z)\ar[r]^{i_*} & \Tt(Y,N'). }$$ 
\end{dfa}

\begin{thma}\label{nondegenere}
Let $X$ and $Y$ be objects of $\Tt$. If there exists a local $\Nn$-cover of $X$ relative to $Y$ and a local $\Nn$-envelope of $\nu X$ relative to $Y$, then the bilienar form $\beta'_{XY}$ constructed in the previous section is non-degenerate.
\end{thma}

For a stronger version of this theorem see also \cite{CR}.

\begin{proof}
Let $f:X\rightarrow Y$ be a morphism in $\Tt$ whose image in $\Tt/\Nn$ is in the kernel of $\beta'$. We have to show that it factorizes through an object of $\Nn$.

Let $p:N\rightarrow X$ be a local $\Nn$-cover of $X$ relative to $Y$, and  let $X'$ be the cone of $p$. The morphism $f$ is in the kernel of $\beta'$, thus for each morphism $g:Y\rightarrow \nu N$ which factorizes through $\nu X'[-1]$, $\beta(fp,g)$ vanishes.
$$\xymatrix{N\ar[r]^p & X\ar[r]\ar[d]^f & X'\ar[r] & N[1] \\ & Y \ar@{-->}[d] \ar[dr]^g & & \\ \nu X[-1] \ar[r] & \nu X'[-1]\ar[r] &\nu N \ar[r] & \nu X }$$
This means that the linear form $\beta(fp,?)$ vanishes on the image of the morphism $\Tt(Y,\nu X'[-1])\longrightarrow \Tt(Y,\nu N). $ This image is canonically isomorphic to the kernel of the morphism
$\Tt(Y,\nu N)\longrightarrow \Tt(Y,\nu X).$

Let $\nu i:\nu X\rightarrow \nu N'$ be a local $\Nn$-envelope of $\nu X$ relative to $Y$. The sequence $$\xymatrix{0\ar[r] & \Tt(Y,\nu X)\ar[r] & \Tt(Y,\nu N')}$$ is then exact. Therefore, the form $\beta(fp,?)$ vanishes on $\Ker (\Tt(Y,\nu N)\longrightarrow \Tt(Y,\nu N')).$
$$\xymatrix{N\ar[r]^p & X\ar[r]\ar[dd]_f \ar[dr]^i & X'\ar[r] & N[1] \\
 &  & N'\ar@{-->}[dl] & \\
 & Y \ar[d]^g &  & \\
  \nu X'[-1]\ar[r] &\nu N \ar[r]\ar@(d,l)[drr] & \nu X\ar[dr]^{\nu i}\ar[r] & \nu X'\\ 
& & &  \nu N' }$$

Now $\beta$ is non-degenerate on $$\Coker (\Tt(N',Y)\longrightarrow \Tt(N,Y))\times \Ker (\Tt(Y,\nu N)\longrightarrow \Tt(Y,\nu N')).$$ Thus the morphism $fp$ lies in $\Coker (\Tt(N',Y)\longrightarrow \Tt(N,Y)),$ that is to say that $fp$ factorizes through $ip$. Since $p:N\rightarrow X$ is a local $\Nn$-cover of $X$,  $f$ factorizes through $N'$.
\end{proof}

\begin{prop}\label{restrictioncouverture}
Let $X$ and $Y$ be objects in $\Tt$. If for each $N$ in $\Nn$ the vector spaces $\Tt(N,X)$ and $\Tt(Y,N)$ are finite-dimensional, then the existence of a local $\Nn$-cover of $X$ relative to $Y$ is equivalent to the existence of a local $\Nn$-envelope of $Y$ relative to $X$.
\end{prop}

\begin{proof}
 Let $g:N\rightarrow X$ be a local $\Nn$-cover of $X$ relative to $Y$. It induces an injection
$$\xymatrix{0\ar[r] & \Tt(X,Y)\ar[r]^{g^*} & \Tt(N,Y).}$$
The space $\Tt(N,Y)$ is finite-dimensional by hypothesis. Fix a basis $(f_1, f_2, \ldots ,f_r)$ of this space. This space is in duality with the space $\Tt(Y,\nu N)$. Let $(f_1',f_2',\ldots ,f_r')$ be the dual basis of the basis $(f_1, f_2, \ldots ,f_r)$. We show that the morphism $$\xymatrix{Y\ar[rr]^(.4){(f_1', \ldots, f_r')} &&\bigoplus_{i=1}^r \nu N}$$ is a local $\Nn$-envelope of $Y$ relative to $X$.
We have a commutative diagram:
$$\xymatrix{\Tt(X,Y)\ar@{>->}[d]^{g^*} \ar[rr]^(.4){(f_1', \ldots, f_r')_*}&&\bigoplus \Tt(X,\nu N)\ar[d]^{g^*}\\ \Tt(N,Y)\ar[rr]^(.4){(f_1', \ldots, f_r')_*} &&\bigoplus \Tt(N, \nu N).}$$ 
If $f$ is in the kernel of $(f_1', \ldots, f_r')_*$, then for all $i=1,\ldots,r$, the morphism $f_i'\circ f\circ g$ is zero. Thus $f\circ g$  is orthogonal on the vectors of the basis $f'_1,\ldots,f'_r$ and therefore vanishes. Since $g$ is a local $\Nn$-cover of $X$ relative to $Y$, $f$ is zero, and the morphism
$$\xymatrix{\Tt(X,Y)\ar[rr]^(.4){(f_1', \ldots, f_r')_*}&&\bigoplus \Tt(X,\nu N)}$$ is injective. Therefore, the morphism $$\xymatrix{Y\ar[rr]^(.4){(f_1', \ldots, f_r')} &&\bigoplus_{i=1}^r \nu N}$$ is a local $\Nn$-envelope of $Y$ relative to $X$. The proof of the converse is dual.
\end{proof}

\subsubsection*{Examples}

 Let $A$ be a finite-dimensional self-injective $k$-algebra. Denote by $\Tt$ the derived category $\Dd^b(\mod A)$
 and by $\Nn$ the triangulated category $\per A$. Since $A$ is finite-dimensional, there is an inclusion $\Nn\subset \Tt$. Moreover $A$ is self-injective so of infinite global dimension. Therefore the inclusion is strict. By \cite{Kel6}, there is an exact sequence of triangulated categories:
 $$\xymatrix{0\ar[r] & \per A\ar[r] & \Dd^b (\mod A)\ar[r] & \underline{\mod} A\ar[r] & 0.}$$
 The derived category $\Dd^b(\mod A)$ admits a Serre functor $\nu=?\lten_A DA$ which stabilizes $\per A$. Thus there is an induced functor in the quotient $\underline{\mod} A$ that we will also denote by $\nu$. Let $\Sigma$ be the suspension of the category $\underline{\mod} A$. One can easily construct (cf. \cite{Ami2}) local $\Nn$-covers and local $\Nn$-envelopes, so we can apply theorem \ref{nondegenere} and   
the functor $\Sigma^{-1}\circ\nu $ is a Serre functor for the stable category $\underline{\mod}A$.

  An article of G. Tabuada \cite{Tab} gives an example of the converse construction. Let $\Cc$ be an algebraic $2$-Calabi-Yau category endowed with a cluster-tilting object. The author constructs a triangulated category $\Tt$ and a triangulated $3$-Calabi-Yau subcategory $\Nn$ such that the quotient category $\Tt/\Nn$ is triangle equivalent to $\Cc$. It is possible to show (cf. \cite{Ami2}) that the categories $\Tt$ and $\Nn$ satisfy the hypotheses of theorem \ref{nondegenere}.

\section{Existence of a cluster-tilting object}

Let $A$ be a differential graded (=dg) $k$-algebra. We denote by $A^e$ the dg algebra $A^{op}\otimes A$.
Suppose that $A$ has the following properties:
\begin{itemize}
 \item $A$ is homologically smooth (\emph{i.e.} the object $A$, viewed as an $A^e$-module, is perfect);
\item for each $p>0$, the space $H^pA$ is zero;
\item the space $H^0A$ is finite-dimensional;
\item $A$ is bimodule $3$-CY, \emph{i.e.} there is an isomorphism in $\Dd(A^e)$  $$R\Hom_{A^e}(A,A^e)\simeq A[-3].$$ 
\end{itemize}

Since $A$ is homologically smooth, the category $\Dd^bA$ is a subcategory of $\per A$ (see \cite{Kel7}, lemma 4.1). We denote by $\pi$ the canonical projection functor $\pi:\per A \rightarrow \Cc=\per A/\Dd^bA$. Moreover, by the same lemma, there is a bifunctorial isomorphism
$$D\Hom_\Dd(L,M)\simeq \Hom_\Dd(M,L[3])$$ for all objects $L$ in $\Dd^b A$ and all objects $M$ in $\per A$. We call this property the \emph{CY property}.

The aim of this section is to show the following result:

\begin{thma}\label{clustertilting}
Let $A$ be a dg $k$-algebra with the above properties. 
  The category $\Cc=\per A/\Dd^b A$ is  $\Hom$-finite and $2$-CY. Moreover, the object $\pi(A)$ is a  cluster-tilting object. Its endomorphism algebra is isomorphic to $H^0A$.
\end{thma}

\subsection{t-structure on $\per A$}

The main tool of the proof of theorem \ref{clustertilting} is the existence of a canonical $t$-structure in $\per A$.

\subsubsection*{$t$-structure on $\Dd A$}
Let $\Dd_{\leq 0}$ be the full subcategory of $\Dd$ whose objects are the dg modules $X$ such that $H^pX$ vanishes for all $p>0$.

\begin{lema}
 The subcategory $\Dd_{\leq 0}$ is an aisle in the sense of Keller-Vossieck \cite{Kel8}.
\end{lema}

\begin{proof}
The canonical morphism $\tau_{\leq 0}A\rightarrow A$ is a quasi-isomorphism of dg algebras. Thus we can assume that $A^p$ is zero for all $p>0$.
 The full subcategory $\Dd_{\leq 0}$ is stable under $X\mapsto X[1]$ and under extensions. We claim that the inclusion $\xymatrix{\Dd_{\leq 0}\ar@{^(->}[r] & \Dd}$ has a right adjoint. Indeed, for each dg A-module $X$, the dg A-module $\tau_{\leq 0}X$ is a dg submodule of $X$, since $A$ is concentrated in negative degrees. Thus $\tau_{\leq 0}$ is a well-defined functor $\Dd\rightarrow \Dd_{\leq 0}$. One can check easily that $\tau_{\leq 0}$ is the right adjoint of the inclusion. 

\end{proof}

\begin{prop}\label{proptstruc}
 Let $\Hh$ be the heart of the $t$-structure, i.e. $\Hh$ is the intersection $\Dd_{\leq 0}\cap\Dd_{\geq 0}$.
We have the following properties:

\hspace{.2in} (i) The functor $H^0$ induces an equivalence from $\Hh$ onto $\Mod H^0A$. 

\hspace{.2in} (ii) For all $X$ and $Y$ in $\Hh$, we have an isomorphism $\Ext^1_{H^0A}(X,Y)\simeq \Hom_\Dd(X,Y[1])$.

\end{prop}
Note that it is not true for general $i$ that $\Ext^i_\Hh(X,Y)\simeq \Hom_\Dd(X,Y[i])$.
\begin{proof}
 (i) We may assume that $A^p=0$ for all $p>0$. We then have a canonical morphism $A\rightarrow H^0A$. The restriction along this morphism yields a functor $\Phi:\Mod H^0A\rightarrow \Hh$ such that $H^0\circ \Phi$ is the identity of $\Mod H^0A$. Thus the functor $H^0:\Hh\rightarrow \Mod H^0A$ is full and essentially surjective. Moreover, it is exact and an object $N\in\Hh$ vanishes if and only if $H^0N$ vanishes. Thus if $f:L\rightarrow M$ is a morphism of $\Hh$ such that $H^0(f)=0$, then $\Imm H^0(f)=0$ implies that $H^0(\Imm f)
=0$ and $\Imm f=0$, so $f=0$. Thus $H^0:\Hh\rightarrow \Mod H^0A$ is also faithful.
 
(ii) By section 3.1.7 of \cite{Bei} there exists a triangle functor $\Dd^b(\Hh)\rightarrow \Dd$ which yields for $X$ and $Y$ in $\Hh$ and for $n\leq1$ an isomorphism (remark (ii) section 3.1.17 p.85) $$\Hom_{\Dd\Hh}(X,Y[n])\simeq \Hom_\Dd(X,Y[n]).$$ 
Applying this for $n=1$ and using (i), we get the result.

\end{proof}

\subsubsection*{$\Hom$-finiteness}

\begin{prop}\label{perAhomfini}
 The category $\per A$ is $\Hom$-finite.
\end{prop}
 
\begin{lema}\label{HpAdimfini}
 For each $p$, the space $H^pA$ is finite-dimensional.
\end{lema}

\begin{proof}
 By hypothesis, $H^pA$ is zero for $p>0$. We prove by induction on $n$ the following statement:
\emph{The space $H^{-n}A$ is finite-dimensional, and for all $p\geq n+1$ the space $\Hom_\Dd(\tau_{\leq -n}A,M[p])$ is finite-dimensional for each $M$ in $\mod H^0A$.}

For $n=0$, the space $H^0A$ is finite-dimensional by hypothesis. Let $M$ be in $\mod H^0A$. Since $\tau_{\leq 0}A$ is ismorphic to $A$, $\Hom_\Dd(\tau_{\leq 0}A,M[p])$ is isomorphic $H^0(M[p])$, and so is zero for $p\geq 1$.

Suppose that the property holds for $n$. Form the triangle:
$$\xymatrix{(H^{-n}A)[n-1]\ar[r] &\tau_{\leq-n-1}A\ar[r] &\tau_{\leq-n}A\ar[r] & (H^{-n}A)[n]}$$
Let $p$ be an integer $\geq n+1$. Applying the functor $\Hom_\Dd(?,M[p])$ we get the long exact sequence:
$$\xymatrix@-1.1pc{\cdots\ar[r] & \Hom_\Dd(\tau_{\leq-n}A,M[p])\ar[r] & \Hom_\Dd(\tau_{\leq -n-1}A,M[p])\ar[r] & \Hom_\Dd((H^{-n}A)[n-1],M[p])\ar[r] & \cdots}.$$

 By induction the space $\Hom_\Dd(\tau_{\leq-n}A,M[p])$ is finite-dimensional. 

Since $M[p]$ is in $\Dd^bA$ we can apply the CY property. If $p$ is $\geq n+3$, we have isomorphisms:
\begin{eqnarray*}
 \Hom_\Dd((H^{-n}A)[n-1],M[p]) &\simeq & \Hom_\Dd((H^{-n}A),M[p-n+1])\\
&\simeq &D\Hom_\Dd(M[p-n-2],H^{-n}A).
\end{eqnarray*}
Since $p-n-2$ is $\geq 1$, this space is zero.

If $p=n+2$, we have the following isomorphisms.
\begin{eqnarray*}
 \Hom_\Dd((H^{-n}A)[n-1],M[n+2]) &\simeq & \Hom_\Dd((H^{-n}A),M[3])\\
&\simeq &D\Hom_\Dd(M,H^{-n}A)\\ & \simeq & D\Hom_{H^0A}(M,H^{-n}A). 
\end{eqnarray*}
The last isomorphism comes from lemma \ref{proptstruc} $(i)$.
By induction, the space $H^{-n}A$ is finite-dimensional. Thus for $p\geq n+2$, the space $\Hom_\Dd((H^{-n}A)[n-1],M[p])$ is finite-dimensional.

If $p=n+1$ we have the following isomorphisms:
\begin{eqnarray*}
 \Hom_\Dd((H^{-n}A)[n-1],M[n+1]) &\simeq & \Hom_\Dd((H^{-n}A),M[2])\\
&\simeq &D\Hom_\Dd(M,H^{-n}A[1])\\ & \simeq & D\Ext^1_{H^0A}(M,H^{-n}A) 
\end{eqnarray*}
The last isomorphism comes from lemma \ref{proptstruc} $(ii)$. By induction, since $H^{-n}A$ is finite-dimensional, the space $\Hom_\Dd((H^{-n}A)[n-1],M[n+1])$ is finite-dimensional and so is $\Hom_\Dd(\tau_{\leq-n-1}A,M[n+1])$.

Now, look at the triangle
$$\xymatrix{\tau_{\leq-n-2}A\ar[r]\ar@{.>}@(d,l)[dr]^0 & \tau_{\leq-n-1}A\ar[r]\ar[d] & (H^{-n-1}A)[n+1]\ar[r]\ar[dl] & (\tau_{\leq-n-2}A)[1]\ar@{.>}@(d,r)[dll]^0\\ & M[n+1]&&}.$$
 The spaces $\Hom_\Dd(\tau_{\leq-n-2}A,M[n+1])$ and $\Hom_\Dd((\tau_{\leq-n-2}A)[1],M[n+1])$ vanish since $M[n+1]$ is in $\Dd_{\geq -n-1}$. Thus we have 
\begin{eqnarray*}
 \Hom_\Dd(\tau_{\leq -n-1}A[n-1],M[n+1]) &\simeq & \Hom_\Dd((H^{-n-1}A)[n+1],M[n+1])\\
&\simeq &\Hom_\Dd(H^{-n-1}A,M)\\ & \simeq & \Hom_{H^0A}(H^{-n-1}A,M).
\end{eqnarray*}
This holds for all finite-dimensional $H^0A$-modules $M$. Thus it holds for the compact cogenerator $M=DH^0A$.  The space $\Hom_{H^0A}(H^{-n-1}A,DH^0A)\simeq DH^{-n-1}A$ is finite-dimensional, and therefore $H^{-(n+1)}A$ is finite-dimensional.
\end{proof}

\begin{proof} (of proposition \ref{perAhomfini})
 For each integer $p$, the space $\Hom_\Dd(A,A[p])\simeq H^p(A)$ is finite-dimensional by lemma \ref{HpAdimfini}. The subcategory of $(\per A)^{op}\times \per A$ whose objects are the pairs $(X,Y)$ such that $\Hom_\Dd(X,Y)$ is finite-dimensional is stable under extensions and passage to direct factors. 
\end{proof}

\subsubsection*{Restriction of the $t$-structure to $\per A$}

\begin{lema}\label{petitlemme}
  For each $X$ in $\per A$, there exist integers $N$ and $M$ such that $X$ belongs to $\Dd_{\leq N}$ and $^\bot\Dd_{\leq M}$.
 \end{lema}
 
 \begin{proof}
 The object $A$ belongs to $\Dd_{\leq 0}$. Moreover, since for $X$ in $\Dd A$, the space $\Hom_\Dd(A,X)$ is isomorphic to $H^0X$, the dg module $A$ belongs to $^\bot \Dd_{\leq-1}$. Thus the property is true for $A$. For the same reasons, it is true for all shifts of $A$.
   Moreover, this property is clearly stable under taking direct summands and extensions. Thus it holds for all objects of $\per A$.
 \end{proof}

 This lemma implies the following result:
 
 \begin{prop}
  The $t$-structure on $\Dd A$ restricts to $\per A$.
 \end{prop}
 
 \begin{proof}
  Let $X$ be in $\per A$, and look at the canonical triangle:
 $$\xymatrix{\tau_{\leq 0} X\ar[r] & X\ar[r] & \tau_{>0} X\ar[r] & (\tau_{\leq 0} X)[1].}$$
Since $\per A$ is $\Hom$-finite, the space $H^pX\simeq \Hom_\Dd(A,X[p])$ is finite-dimensional for all $p\in\mathbb{Z}$ and vanishes for all $p\gg 0$ by lemma \ref{petitlemme}. Thus the object $\tau_{>0}X$ is in $\Dd^bA$ and so in $\per A$. Since $\per A$ is a triangulated subcategory, it follows that $\tau_{\leq 0}X$ also lies in $\per A$. 
 \end{proof}

\begin{prop}
 Let $\pi$ be the projection $\pi:\per A\rightarrow \Cc$. Then for $X$ and $Y$ in $\per A$, we have
 $$\Hom_\Cc(\pi X,\pi Y)=\lim_{\rightarrow}\Hom_\Dd(\tau_{\leq n}X,\tau_{\leq n}Y)$$   
\end{prop}
\begin{proof} Let $X$ and $Y$ be in $\per A$. An element of $\underset{\rightarrow}{\lim}\Hom_\Dd(\tau_{\leq n}X,\tau_{\leq n}Y)$ is an equivalence class of morphisms $\tau_{\leq n}X\rightarrow \tau_{\leq n}Y$. Two morphisms $f:\tau_{\leq n}X\rightarrow \tau_{\leq n}Y$ and $g:\tau_{\leq m}X\rightarrow \tau_{\leq m}Y$ with $m\geq n$ are equivalent if there is a commutative square:
$$\xymatrix{\tau_{\leq n}X\ar[r]^f\ar[d] & \tau_{\leq n}Y\ar[d]\\ \tau_{\leq m}X\ar[r]^g & \tau_{\leq m}Y}$$
 where the vertical arrows are the canonical morphisms.
 If $f$ is a morphism $f:\tau_{\leq n}X\rightarrow \tau_{\leq n}Y$, we can form the following morphism from $X$ to $Y$ in $\Cc$:
 $$\xymatrix{&\tau_{\leq n}X\ar[dl]\ar@{.>}[r]^f\ar[dr]& \tau_{\leq n}Y\ar@{.>}[d]\\ X&&Y,}$$ where the morphisms $\tau_{\leq n}X\rightarrow X$ and $\tau_{\leq n} Y\rightarrow Y$ are the canonical morphisms. This is a morphism from $\pi X$ to $\pi Y$ in $\Cc$ because the cone of the morphism $\tau_{\leq n}X\rightarrow X$ is in $\Dd^bA$. Moreover, if $f:\tau_{\leq n}X\rightarrow \tau_{\leq n}Y$ and $g:\tau_{\leq m}X\rightarrow \tau_{\leq m}Y$ are equivalent, there is an equivalence of diagrams:
 $$\xymatrix{&\tau_{\leq n}X\ar[dd]\ar[dl]\ar@{.>}[r]^f\ar[dr]& \tau_{\leq n}Y\ar@{.>}[d]\ar@(r,r)@{.>}[dd]\\ X&&Y\\&\tau_{\leq m}X\ar[ul]\ar@{.>}[r]^g\ar[ur]& \tau_{\leq m}Y\ar@{.>}[u]}$$ Thus we have a well-defined map from $\underset{\rightarrow}{\lim}\Hom_\Dd(\tau_{\leq n}X,\tau_{\leq n}Y)$ to $\Hom_\Cc(\pi X, \pi Y)$ which is injective. 
 
 Now let $\xymatrix@-1.5pc{&X'\ar[dr]^s\ar[dl]&\\X&&Y}$ be a morphism in $\Hom_\Cc(\pi X,\pi Y)$. Let $X''$ be the cone of $s$. It is an object of $\Dd^bA$, and therefore lies in $\Dd_{>n}$ for some $n\ll 0$. Thus there are no morphisms from $\tau_{\leq n}X$ to $X''$ and there is a factorization:
 $$\xymatrix{& \tau_{\leq n}X\ar[d]\ar[dr]^0\ar@{.>}[dl] & &\\X'\ar[r]^s&X\ar[r]&X''\ar[r] & X'[1]}$$
 We obtain  an isomorphism of diagrams:
 $$\xymatrix@-1.3pc{&X'\ar[dr]\ar[dl]_s&\\X&&Y\\&\tau_{\leq n}X\ar[ur]_f\ar[ul]\ar[uu]}$$
 The morphism $f:\tau_{\leq n}X\rightarrow Y$ induces a morphism $f':\tau_{\leq n}X\rightarrow \tau_{\leq n} Y$ which lifts the given morphism. Thus the map from $\underset{\rightarrow}{\lim}\Hom_\Dd(\tau_{\leq n}X,\tau_{\leq n}Y)$ to $\Hom_\Cc(\pi X, \pi Y)$ is surjective.
\end{proof}

\subsection{Fundamental domain}

Let $\Ff$ be the following subcategory of $\per A$:
$$\Ff=\Dd_{\leq 0}\cap \;^\bot\! \Dd_{\leq -2}\cap \per A.$$
The aim of this section is to show:
\begin{prop}\label{equivalenceFC}
 The projection functor $\pi:\per A \rightarrow \Cc$ induces a $k$-linear equivalence between $\Ff$ and $\Cc$. 
\end{prop}

\subsubsection*{$add(A)$-approximation for objects of the fundamental domain}
\begin{lema}\label{addAapprox}
 For each object $X$ of $\Ff$, there exists a triangle $$\xymatrix{P_1\ar[r] &P_0\ar[r] &X\ar[r] & P_1[1]}$$ with $P_0$ and $P_1$ in $add(A)$.
\end{lema}

\begin{proof}
 For $X$ in $\per A$, the morphism $$\begin{array}{rcl}\Hom_\Dd(A,X)&\rightarrow &\Hom_\Hh(H^0A,H^0X)\\f&\mapsto & H^0(f)\end{array} $$ is an isomorphism since $\Hom_\Dd(A,X)\simeq H^0X$. 
Thus it is possible to find a morphism $P_0\rightarrow X$, with $P_0$ a free dg $A$-module, inducing an epimorphism $\xymatrix{H^0P_0\ar@{->>}[r] &H^0X}$. Now take $X$ in $\Ff$ and $P_0\rightarrow X$ as previously and form the triangle $$\xymatrix{P_1\ar[r] & P_0\ar[r] & X\ar[r] & P_1[1].}$$

\textit{Step 1: The object $P_1$ is in $\Dd_{\leq 0}\cap\; ^\bot\! \Dd_{\leq-1}$.}\smallskip\\
The objects $X$ and $P_0$ are in $\Dd_{\leq 0}$, so $P_1$ is in $\Dd_{\leq 1}$. Moreover, since $H^0P_0\rightarrow H^0X$ is an epimorphism, $H^1(P_1)$ vanishes and $P_1$ is in $\Dd_{\leq 0}$.

Let $Y$ be in $\Dd_{\leq -1}$, and look at the long exact sequence:
$$\xymatrix{\cdots\ar[r]  & \Hom_\Dd(P_0,Y)\ar[r] & \Hom_\Dd( P_1,Y)\ar[r] & \Hom_\Dd(X[-1],Y)\ar[r] & \cdots}.$$
The space $\Hom_\Dd(X[-1],Y)$ vanishes since $X$ is in $^\bot \Dd_{\leq -2}$ and $Y$ is in $\Dd_{\leq -1}$. The object $P_0$ is free, and $H^0Y$ is zero, so the space $\Hom_\Dd(P_0,Y)$ also vanishes. Consequently, the object $P_1$ is in $^\bot\Dd_{\leq -1}$.\medskip

\textit{Step 2: $H^0P_1$ is a projective $H^0A$-module.}\smallskip\\
Since $P_1$ is in $\Dd_{\leq 0}$ there is a triangle
$$\xymatrix{\tau_{\leq -1}P_1\ar[r] &P_1\ar[r] & H^0P_1\ar[r] & (\tau_{\leq -1}P_1)[1].}$$
Now take an object $M$ in the heart $\Hh$, and look at the long exact sequence:
$$\xymatrix{\cdots\ar[r]  & \Hom_\Dd((\tau_{\leq -1}P_1)[1],M[1])\ar[r] & \Hom_\Dd( H^0P_1,M[1])\ar[r] & \Hom_\Dd(P_1,M[1])\ar[r] & \cdots}.$$
The space $\Hom_\Dd((\tau_{\leq -1}P_1)[1],M[1])$ is zero because $\Hom_\Dd(\Dd_{\leq-2},\Dd_{\geq -1})$ vanishes in a $t$-structure. Moreover, the space $\Hom_\Dd(P_1,M[1])$ vanishes because $P_1$ is in $^\bot\Dd_{\leq -1}$.
Thus $\Hom_\Dd( H^0P_1,M[1])$ is zero. But this space is isomorphic to the space $\Ext^1_\Hh(H^0P_1,M)$ by proposition \ref{proptstruc}. This proves that $H^0P_1$ is a projective $H^0A$-module.\medskip
 
\textit{Step 3: $P_1$ is isomorphic to an object of $add(A)$.}\smallskip\\
As previously, since $H^0P_1$ is projective, it is possible to find an object $P$ in $add(A)$ and a morphism $P\rightarrow P_1$ inducing an isomorphism $H^0P\rightarrow H^0P_1 $. Form the triangle 
$$\xymatrix{Q\ar[r]^u & P\ar[r]^v & P_1\ar[r]^w & Q[1]}$$
Since $P$ and $P_1$ are in $\Dd_{\leq 0}$ and $H^0(v)$ is surjective, the cone $Q$ lies in $\Dd_{\leq 0}$. But then $w$ is zero since $P_1$ is in $^\bot \Dd_{\leq -1}$. Thus the triangle splits, and $P$ is isomorphic to the direct sum $P_1\oplus Q$. Therefore we have a short exact sequence:
$$\xymatrix{0\ar[r] & H^0Q\ar[r] & H^0P\ar[r] & H^0P_1\ar[r] & 0,}$$
and $H^0Q$ vanishes. The object $Q$ is in $\Dd_{\leq-1}$, the triangle splits, and there is no morphism between $P$ and $\Dd_{\leq -1}$, so $Q$ is the zero object.  

\end{proof}

\subsubsection*{Equivalence between the shifts of $\Ff$}
\begin{lema}\label{equivalenceF}
 The functor $\tau_{\leq -1}$ induces an equivalence from $\Ff$ to $\Ff[1]$
\end{lema}

\begin{proof}

\textit{Step 1: The image of the functor $\tau_{\leq-1}$ restricted to $\Ff$ is in $\Ff[1]$.}\smallskip\\
Recall that $\Ff$ is $\Dd_{\leq 0}\cap \;^\bot\! \Dd_{\leq -2}\cap \per A$ so $\Ff[1]$ is $\Dd_{\leq -1}\cap \;^\bot\! \Dd_{\leq -3}\cap \per A$.
 Let $X$ be an object in $\Ff$. By definition, $\tau_{\leq -1}X$ lies in $\Dd_{\leq-1}$ and there is a canonical triangle:
$$\xymatrix{\tau_{\leq-1}X\ar[r]& X\ar[r] & H^0X\ar[r] & \tau_{\leq -1}X[1]}.$$
Now let $Y$ be an object in $\Dd_{\leq -3}$ and form the long exact sequence
$$\xymatrix{\cdots\ar[r]  & \Hom_\Dd(X,Y)\ar[r] & \Hom_\Dd( \tau_{\leq-1}X,Y)\ar[r] & \Hom_\Dd((H^0X)[-1],Y)\ar[r] & \cdots}$$
Since $X$ is in $^\bot\Dd_{\leq-2}$, the space $\Hom_\Dd(X,Y)$ vanishes. The object $H^0X[-1]$ is of finite total dimension, so by the CY property, we have an isomorphism
$$\Hom_\Dd(H^0X[-1],Y)\simeq D\Hom_\Dd(Y,H^0X[2]).$$ But since $\Hom_\Dd(\Dd_{\leq-3},\Dd_{\geq-2})$ is zero, the space $\Hom_\Dd((H^0X)[-1],Y)$ vanishes and $\tau_{\leq-1}X$ lies in $^\bot\Dd_{\leq -3}$.\medskip

\textit{Step 2: The functor $\tau_{\leq -1}:\Ff\rightarrow \Ff[1]$ is fully faithful.}\smallskip\\
Let $X$ and $Y$ be two objects in $\Ff$ and $f:\tau_{\leq-1}X\rightarrow \tau_{\leq-1}Y$ be a morphism.
$$\xymatrix{H^0X[-1]\ar[r] & \tau_{\leq-1}X\ar[r]\ar[d]^f & X\ar@{.>}[d]\ar[r] & H^0X \\  H^0Y[-1]\ar[r] & \tau_{\leq-1}Y\ar[r]^i & Y\ar[r] & H^0Y}$$
The space $\Hom_\Dd(H^0X[-1],Y)$ is isomorphic to $D\Hom_\Dd(Y,H^0X[2])$ by the CY property. Since $Y$ is in $^\bot \Dd_{\leq -2}$, this space is zero, and  the composition $i\circ f$ factorizes through the canonical morphism $\tau_{\leq -1}X\rightarrow X$. Therefore, the functor $\tau_{\leq-1}$ is full.

Let $X$ and $Y$ be objects of $\Ff$ and $f:X\rightarrow Y$ a morphism satisfying $\tau_{\leq-1}f=0$. It induces a morphism of triangles:
$$\xymatrix{H^0X[-1]\ar[r] \ar[d]& \tau_{\leq-1}X\ar[r]^i\ar[d]^0& X\ar[r] \ar[d]^f& H^0X\ar@{.>}[ld]\ar[d] \\  H^0Y[-1]\ar[r] & \tau_{\leq-1}Y\ar[r]& Y\ar[r] & H^0Y}$$
The composition $f\circ i$ vanishes, so $f$ factorizes through $H^0X$. But by the CY property the space of morphisms $\Hom_\Dd(H^0X,Y)$ is isomorphic to $D\Hom_\Dd(Y,H^0X[3])$ which is zero since $Y$ is in $^\bot\Dd_{\leq -2}$. Thus the functor $\tau_{\leq-1}$ restricted to $\Ff$ is faithful.\medskip

\textit{Step 3: The functor $\tau_{\leq -1}:\Ff\rightarrow \Ff[1]$ is essentially surjective.}\smallskip\\
Let $X$ be in $\Ff[1]$. By the previous lemma there exists a triangle
$$\xymatrix{P_1[1]\ar[r] & P_0[1]\ar[r] & X\ar[r] & P_1[2]}$$
with $P_0$ and $P_1$ in $add(A)$. Denote by $\nu$ the Nakayama functor on the projectives of $\mod H^0A$. Let $M$ be the kernel of the morphism $\nu H^0P_1\rightarrow \nu H^0 P_0$. It lies in the heart $\Hh$.\medskip

\textit{Substep (i): There is an isomorphism of functors: $\Hom(?,X[1])_{|_\Hh}\simeq\Hom_\Hh(?,M)$}\smallskip\\
Let $L$ be in $\Hh$. We then have a long exact sequence:
$$\xymatrix@-1.1pc{\cdots\ar[r]  & \Hom_\Dd(L,P_0[2])\ar[r] & \Hom_\Dd( L,X[1])\ar[r] & \Hom_\Dd(L,P_1[3])\ar[r] &\Hom_\Dd(L,P_0[3])\ar[r] & \cdots}.$$ 
The space $\Hom_\Dd(L,P_0[2])$ is isomorphic to $D\Hom_\Dd(P_0,L[1])$ by the CY property, and vanishes because $P_0$ is in $^\bot\Dd_{\leq-1}$.
Moreover, we have the following isomorphisms:
\begin{eqnarray*}
 \Hom_\Dd(L,P_1[3]) & \simeq & D\Hom_\Dd(P_1,L)\\ & \simeq & D\Hom_\Hh(H^0P_1,L)\\ & \simeq & \Hom_\Hh(L,\nu H^0P_1).
\end{eqnarray*}
Thus $\Hom_\Dd(?,X[1])_{|_\Hh}$ is isomorphic to the kernel of $\Hom_\Hh(?,\nu H^0P_1)\rightarrow \Hom_\Hh(?,\nu H^0P_0)$, which is just $\Hom_\Hh(?,M)$.\medskip

\textit{Substep (ii): There is a monomorphism of functors:$\xymatrix@-1.3pc{\Ext^1_\Hh(?,M)\ar@{^(->}[r] & \Hom_\Dd(?,X[2])_{|_\Hh}}$.}\smallskip\\ 
For $L$ in $\Hh$, look at the following long exact sequence:
$$\xymatrix@-1.1pc{\cdots\ar[r]  & \Hom_\Dd(L,P_1[3])\ar[r] & \Hom_\Dd( L,P_1[3])\ar[r] & \Hom_\Dd(L,X[2])\ar[r] &\Hom_\Dd(L,P_1[4])\ar[r] & \cdots}.$$ 
The space $\Hom_\Dd(L,P_1[4])$ is isomorphic to $D\Hom_\Dd(P_1[1],L)$ which is zero since $P_1[1]$ is in $\Dd_{\leq-1}$ and $L$ is in $\Dd_{\geq 0}$. Thus the functor $\Hom_\Dd(?,X[2])_{|_\Hh}$ is isomorphic to the cokernel of $\Hom_\Hh(?,\nu H^0P_1)\rightarrow \Hom_\Hh(?,\nu H^0P_0)$.
By defninition, $\Ext^1_\Hh(?,M)$ is the first homology of a complex of the form:
$$\xymatrix{\cdots\ar[r] & 0\ar[r] &\Hom_\Hh(?,\nu H^0P_1)\ar[r] & \Hom_\Hh(?,\nu H^0P_0)\ar[r] & \Hom_\Hh(?,I)\ar[r] &\cdots},$$ where $I$ is an injective $H^0A$-module.
Thus we get the canonical injection:
$$\xymatrix{\Ext^1_\Hh(?,M)\ar@{^(->}[r] & \Hom_\Dd(?,X[2])_{|_\Hh}.}$$

Now form the following triangle:
$$\xymatrix{X\ar[r] & Y\ar[r] & M\ar[r] & X[1].}$$

\textit{Substep (iii): $Y$ is in $\Ff$ and $\tau_{\leq-1}Y$ is isomorphic to $X$.}\smallskip\\
Since $X$ and $M$ are in $\Dd_{\leq 0}$, $Y$ is in $\Dd_{\leq 0}$. Let $Z$ be in $\Dd_{\leq -2}$ and form the following long exact sequence:
$$\xymatrix@-1.1pc{\cdots \Hom_\Dd(X[1],Z)\ar[r] & \Hom_\Dd(M,Z)\ar[r] & \Hom_\Dd(Y,Z)\ar[r] &\Hom_\Dd(X,Z)\ar[r] &\Hom_\Dd(M[-1],Z)\cdots.}$$
By the CY property and the fact that $Z[2]$ is in $\Dd_{\leq 0}$, we have isomorphisms
 \begin{eqnarray*}
 \Hom_\Dd(M[-1],Z) & \simeq & D\Hom_\Dd(Z[-2],M)\\ & \simeq & D\Hom_\Hh(H^{-2}Z,M).
\end{eqnarray*}
Moreover, since $X$ is in $^\bot\Dd_{\leq-3}$, we have
\begin{eqnarray*}
 \Hom_\Dd(X,Z) & \simeq & \Hom_\Dd(X,(H^{-2}Z)[2])\\ & \simeq & D\Hom_\Hh(H^{-2}Z,X[1]).
\end{eqnarray*}
By substep (i) the functors $\Hom_\Hh(?,M)$ and $\Hom_\Dd(?,X[1])_{|_\Hh}$ are isomorphic. Therefore we deduce that the morphism $\Hom_\Dd(X,Z)\rightarrow\Hom_\Dd(M[-1],Z)$ is an isomorphism. 

Now look at the triangle 
$$\xymatrix{\tau_{\leq-3}Z\ar[r] & Z\ar[r] & H^{-2}Z[2]\ar[r] & (\tau_{\leq-3}Z)[1]}$$
and form the commutative diagram
$$\xymatrix@C=.7pc{ \Hom_\Dd(M,\tau_{\leq-3}Z)\ar[r] & \Hom_\Dd(M,Z)\ar[r] & \Hom_\Dd(M,H^{-2}Z[2])\ar[r] &\Hom_\Dd(M,\tau_{\leq-3}Z[1]) \\
 \Hom_\Dd(X[1],\tau_{\leq-3}Z)\ar[r]\ar[u]^a & \Hom_\Dd(X[1],Z)\ar[r]\ar[u]^b & \Hom_\Dd(X[1],H^{-2}Z[2])\ar[r]\ar[u]^c &\Hom_\Dd(X[1],\tau_{\leq-3}Z[1])\ar[u]^d}.$$
By the CY property and the fact that $(\tau_{\leq-3}Z)[-3]$ is in $\Dd_{\leq 0}$, we have isomorphisms
 \begin{eqnarray*}
 \Hom_\Dd(M[-1],\tau_{\leq -3}Z[-1]) & \simeq & D\Hom_\Dd(\tau_{\leq-3}Z[-3],M)\\ & \simeq & D\Hom_\Hh(H^{-3}Z,M).
\end{eqnarray*}
Since $X$ is in $^\bot\Dd_{\leq-3}$, we have
\begin{eqnarray*}
 \Hom_\Dd(X,(\tau_{\leq-3}Z)[-1]) & \simeq & \Hom_\Dd(X,H^{-3}Z[-2])\\ & \simeq & D\Hom_\Hh(H^{-3}Z,X[1]).
\end{eqnarray*}
Now we deduce from substep (i) that $a[-1]$ is an isomorphism.

 The space $\Hom_\Dd(X[1],\tau_{\leq-3}Z[1])$ is zero because $X$ is $^\bot\Dd_{\leq-3}$. Moreover there are isomorphisms
 \begin{eqnarray*}
 \Hom_\Dd(M,H^{-2}Z[2]) & \simeq & D\Hom_\Dd(H^{-2}Z,M[1])\\ & \simeq & D\Ext^1_\Hh(H^{-2}Z,M).
\end{eqnarray*}
The space $\Hom_\Dd(X[1],H^{-2}Z[2])$ is isomorphic to $D\Hom_\Dd(H^{-2}Z,X[2])$. And by substep (ii), the morphism $\Ext^1_\Hh(?,M)\rightarrow\Hom_\Dd(?,X[2])_{|_\Hh}$ is injective, so $c$ is surjective.
Therefore using a weak form of the five-lemma we deduce that $b$ is surjective.

 Finally, we have the following exact sequence:
$$\xymatrix@-.9pc{ \Hom_\Dd(X[1],Z)\ar@{->>}[r] & \Hom_\Dd(M,Z)\ar[r] & \Hom_\Dd(Y,Z)\ar[r] &\Hom_\Dd(X,Z)\ar[r]^(.4){\sim} &\Hom_\Dd(M[-1],Z)}$$
Thus the space $\Hom_\Dd(M,Z)$ is zero, and $Z$ is in $^\bot\Dd_{\leq-2}$.

It is now easy to see that there is an isomorphism of triangles:
$$\xymatrix{\tau_{\leq-1}Y\ar[r]\ar[d]& Y\ar[r]\ar@{=}[d] & H^0Y\ar[r]\ar[d] &\tau_{\leq-1}Y[1]\ar[d]\\X\ar[r] & Y\ar[r] & M\ar[r] & X[1].}$$

\end{proof}

\subsubsection*{Proof of proposition \ref{equivalenceFC}}

\textit{Step 1: The functor $\pi$ restricted to $\Ff$ is fully faithful.}\smallskip\\
Let $X$ and $Y$ be objects in $\Ff$. By proposition \ref{proptstruc} (iii), the space $\Hom_\Cc(\pi X,\pi Y)$ is isomorphic to the direct limit $\underset{\rightarrow}{\lim}\Hom_\Dd(\tau_{\leq n}X,\tau_{\leq n} Y)$. A morphism between $X$ and $Y$ in $\Cc$ is a diagram of the form 
$$\xymatrix@-1.5pc{&\tau_{\leq n}X\ar[dr]\ar[dl]&\\ X&&Y.}$$ The canonical triangle $$\xymatrix{(\tau_{>n}X)[-1]\ar[r] & \tau_{\leq n}X\ar[r] & X\ar[r] & \tau_{>n}X}$$ yields a long exact sequence:
$$\xymatrix@-1.2pc{\cdots\ar[r] &\Hom_\Dd(\tau_{>n}X,Y)\ar[r] & \Hom_\Dd(X,Y)\ar[r] & \Hom_\Dd(\tau_{\leq n}X,Y)\ar[r] & \Hom_\Dd((\tau_{> n}X)[-1],Y)\ar[r] & \cdots}$$
The space $\Hom_\Dd((\tau_{>n}X)[-1],Y)$ is isomorphic to the space $D\Hom_\Dd(Y,(\tau_{>n}X)[2])$. The object $X$ is in $\Dd_{\leq 0}$, thus so is $\tau_{>n}X$, and the space $D\Hom_\Dd(Y,(\tau_{>n}X)[2])$ vanishes. For the same reasons, the space $\Hom_\Dd(\tau_{>n}X,Y)$ vanishes. Thus there are bijections 
$$\xymatrix{ \Hom_\Dd(\tau_{\leq n}X,\tau_{\leq n}Y)\ar[r]^\sim &\Hom_\Dd(\tau_{\leq n}X,Y)\ar[r]^\sim & \Hom_\Dd(X,Y)}$$
Therefore, the functor $\pi:\Ff\rightarrow \Cc$ is fully faithful.\medskip

\textit{Step 2: For $X$ in $\per A$, there exists an integer $N$ and an object $Y$ of $\Ff[-N]$ such that $\pi X$ and $\pi Y$ are isomorphic in $\Cc$.}\smallskip\\
Let $X$ be in $\per A$. By lemma \ref{petitlemme}, there exists an integer $N$ such that $X$ is in $^\bot\Dd_{\leq N-2}$. For an object $Y$ in $\Dd_{\leq N-2}$, the space $\Hom_\Dd((\tau_{> N}X)[-1],Y)$ is isomorphic to $D\Hom_\Dd(Y,(\tau_{> N}X)[2])$ and thus vanishes. Therefore, $\tau_{\leq N}X$ is still in $^\bot\Dd_{\leq N-2}$, and thus is in $\Ff[-N]$. Since $\tau_{> N}X$ is in $\Dd^bA$, the objects $\tau_{\leq N}X$ and $X$ are isomorphic in $\Cc$. \medskip

\textit{Step 3: The functor $\pi$ restricted to $\Ff$ is essentially surjective.}\smallskip\\
Let $X$ be in $\per A$ and $N$ such that $\tau_{\leq N}X$ is in $\Ff[-N]$. By lemma \ref{equivalenceF}, $\tau_{\leq-1}$ induces an equivalence between $\Ff$ and $\Ff[1]$. Thus since the functor $\pi\circ\tau_{\leq-1}:\per A\rightarrow \Cc$ is isomorphic to $\pi$, there exists an object $Y$ in $\Ff$ such that $\pi(Y)$ and $\pi(X)$ are isomorphic in $\Cc$. Therefore, the functor $\pi$ restricted to $\Ff$ is essentially surjective.

\begin{prop}\label{suiteexacte}
 If $X$ and $Y$ are objects in $\Ff$, there is a short exact sequence:
$$\xymatrix{0\ar[r] & \Ext^1_\Dd(X,Y)\ar[r] & \Ext^1_\Cc(X,Y)\ar[r] & D\Ext^1_\Dd(Y,X)\ar[r] & 0.}$$
\end{prop}

\begin{proof}
Let $X$ and $Y$ be in $\Ff$. The canonical triangle
$$\xymatrix{\tau_{<0}X\ar[r] & X\ar[r]& \tau_{\geq 0}X\ar[r] & (\tau_{<0}X)[1]}$$ yields the long exact sequence:
$$\xymatrix@-1.1pc{ \Hom_\Dd((\tau_{\geq 0}X)[-1],Y[1])& \Hom_\Dd(\tau_{<0}X,Y[1])\ar[l] & \Hom_\Dd(X,Y[1])\ar[l] & \Hom_\Dd(\tau_{\geq 0}X,Y[1])\ar[l]}.$$

The space $\Hom_\Dd(X[-1],Y[1])$ is zero because $X$ is in $^\bot\Dd_{\leq-2}$ and $Y$ is in $\Dd_{\leq 0}$. Moreover, the space $\Hom_\Dd(\tau_{\geq 0}X,Y[1])$ is zero because of the CY property. Thus this long sequence reduces to a short exact sequence:
$$\xymatrix{0\ar[r] & \Ext^1_\Dd(X,Y)\ar[r] & \Hom_\Dd(\tau_{<0}X,Y[1])\ar[r]& \Hom_\Dd((\tau_{\geq 0}X)[-1],Y[1])\ar[r] & 0}.$$
\textit{Step 1: There is an isomorphism $\Hom_\Dd((\tau_{\geq 0}X)[-1],Y)\simeq D\Ext^1_\Dd(Y,X).$}\smallskip\\
The space $\Hom_\Dd((\tau_{\geq 0}X)[-1],Y[1])$ is isomorphic to $D\Hom_\Dd(Y,\tau_{\geq 0}X[1])$ by the CY property.
$$\xymatrix{&Y\ar@{..>}[dl]_0\ar[d]\ar[dr]\ar@(r,u)@{..>}[drr]^0&&\\(\tau_{<0}X)[1]\ar[r] & X[1]\ar[r]& (\tau_{\geq 0}X)[1]\ar[r] & (\tau_{<0}X)[2]}$$
  But since $\Hom_\Dd(Y,(\tau_{<0}X)[1])$ and $\Hom_\Dd(Y,(\tau_{<0}X)[2])$ are zero, we have an isomorphism $$\Hom_\Dd(\tau_{\geq 0}X[-1],Y)\simeq D\Ext^1_\Dd(Y,X).$$
\textit{Step 2: There is an isomorphism $\Ext^1_\Cc(\pi X,\pi Y)\simeq\Hom_\Dd(\tau_{\leq-1}X,Y[1])$.}\smallskip\\
By lemma \ref{equivalenceF}, the object $\tau_{<0}X$ belongs to $\Ff[1]$ and clearly $Y[1]$ belongs to $\Ff[1]$. By proposition \ref{equivalenceFC} (applied to the shifted $t$-structure), the functor $\pi:\per A\rightarrow \Cc$ induces an equivalence from $\Ff[1]$ to $\Cc$ and clearly we have $\xymatrix{\pi(\tau_{<0}X,Y[1])\ar[r]^\sim & \pi(X).}$ We obtain bijections
$$\xymatrix{\Hom_\Dd(\tau_{<0}X,Y[1])\ar[r]^\sim & \Hom_\Dd(\pi\tau_{<0}X,\pi Y[1])\ar[r]^\sim&\Hom_\Dd(\pi X,\pi Y[1]).}$$  
\end{proof}

\subsubsection*{Proof of theorem \ref{clustertilting}}

\textit{Step 1: The category $\Cc$ is $\Hom$-finite and $2$-CY.}\smallskip\\
The category $\Ff$ is obviously $\Hom$-finite, hence so is $\Cc$ by proposition \ref{equivalenceFC}.
The categories $\Tt=\per A$ and $\Nn=\Dd^bA\subset \per A$ satisfy the hypotheses of section 1. By \cite{Kel7}, thanks to the CY property,  there is a bifunctorial non degenerate bilinear form:
$$\beta_{N,X}:\Hom_\Dd(N,X)\times \Hom_\Dd(X,N[3])\rightarrow k$$
for $N$ in $\Dd^bA$ and  $X$ in $\per A$.
Thus, by section 1, there exists a bilinear bifunctorial form 
$$\beta'_{X,Y}:\Hom_\Cc(X,Y)\times \Hom_\Cc(Y, X[2])\rightarrow k$$
for $X$ and $Y$ in $\Cc=\per A/\Dd^bA$. We would like to show that it is non degenerate. Since $\per A$ is $\Hom$-finite, by theorem \ref{nondegenere} and proposition \ref{restrictioncouverture}, it is sufficient to show the existence of local $\Nn$-envelopes. Let $X$ and $Y$ be objects of $\per A$. Therefore by lemma \ref{petitlemme}, $X$ is in $^\bot\Dd_{\leq N}$. Thus there is an injection
$$\xymatrix{0\ar[r]&\Hom_\Dd(X,Y)\ar[r]&\Hom_\Dd(X,\tau_{>N}Y)}$$ 
and $Y\rightarrow \tau_{>N}Y$ is a local $\Nn$-envelope relative to $X$.  Therefore, $\Cc$ is $2$-CY.

\medskip

\textit{Step 2: The object $\pi A$ is a cluster-tilting object of the category $\Cc$.}\smallskip\\
Let $A$ be the free dg $A$-module in $\per A$. Since $H^1A$ is zero, the space $\Ext^1_\Dd(A,A)$ is also zero. Thus by the short exact sequence 
$$\xymatrix{0\ar[r] & \Ext^1_\Dd(A,A)\ar[r] & \Ext^1_\Cc(\pi A,\pi A)\ar[r]& D\Ext^1_\Dd(A,A)\ar[r]& 0}$$ of
 proposition \ref{suiteexacte}, $\pi(A)$ is a rigid object of $\Cc$. 
Now let $X$ be an object of $\Cc$. By proposition \ref{equivalenceFC}, there exists an object $Y$ in $\Ff$ such that $\pi Y$ is isomorphic to $X$. Now by lemma \ref{addAapprox} , there exists a triangle in $\per A$ $$\xymatrix{P_1\ar[r]&P_0\ar[r] & Y\ar[r] & P_1[1]}$$ with $P_1$ and $P_0$ in $add(A)$. Applying the triangle functor $\pi$  we get a triangle in $\Cc$:
$$\xymatrix{\pi P_1\ar[r] & \pi P_0\ar[r] &  X \ar[r]& \pi P_1[1]}$$ with $\pi P_1$ and $\pi P_0$ in $add(\pi A)$. If $\Ext^1_\Cc(\pi A,X)$ vanishes, this triangle splits and $X$ is a direct factor of $\pi P_0$. Thus, the object $\pi A$ is a cluster-tilting object in the 2-CY category $\Cc$.

\section{Cluster categories for Jacobi-finite quivers with potential}

\subsection{Ginzburg dg algebra}

Let  $Q$ be a finite quiver. For each arrow $a$ of $Q$, we define the \emph{cyclic derivative with respect to} $a$ $\partial_a$ as the unique linear map
$$\partial_a:kQ/[kQ,kQ]\rightarrow kQ$$
which takes the class of a path $p$ to the sum $\sum_{p=uav}vu $ taken over all decompositions of the path $p$ (where $u$ and $v$ are possibly idempotents $e_i$ associated to a vertex $i$ of $Q$). 

An element $W$ of $kQ/[kQ,kQ]$ is called a \emph{potential on }$Q$. It is given by  a linear combination of cycles in $Q$. 

\begin{dfa}[Ginzburg]\cite{Gin}(section 4.2)
 Let $Q$ be a finite quiver and $W$ a potential on $Q$. Let $\widehat{Q}$ be the graded quiver with the same vertices as $Q$ and whose arrows are
\begin{itemize}
 \item the arrows of $Q$ (of degree $0$),
\item an arrow $a^*:j\rightarrow i$ of degree $-1$ for each arrow $a:i\rightarrow j$ of $Q$,
\item a loop $t_i:i\rightarrow i$ of degree $-2$ for each vertex $i$ of $Q$.
\end{itemize}
The \emph{Ginzburg dg algebra} $\Gamma(Q,W)$ is a dg $k$-algebra whose underlying graded algebra is the graded path algebra $k\widehat{Q}$. Its differential is the unique linear endomorphism homogeneous of degree $1$ which satisfies  the Leibniz rule
$$d(uv)=(du)v+(-1)^pudv,$$
for all homogeneous $u$ of degree $p$ and all $v$, and takes the following values on the arrows of $\widehat{Q}$:
\begin{itemize}
 \item $da=0$ for each arrow $a$ of $Q$,
\item $d(a^*)=\partial_aW$ for each arrow $a$ of $Q$,
\item $d(t_i)=e_i(\sum_{a}[a,a^*])e_i$ for each vertex $i$ of $Q$ where $e_i$ is the idempotent associated to $i$ and the sum runs over all arrows of $Q$.
\end{itemize}

\end{dfa}

The strictly positive homology of this dg algebra clearly vanishes. Moreover B. Keller showed the following result:

\begin{thma}[Keller]\cite{Kel10}\label{ginzburg3CY}
 Let $Q$ be a finite quiver and $W$ a potential on $Q$. Then the Ginzburg dg algebra $\Gamma(Q,W)$ is  homologically smooth and bimodule $3$-CY.
\end{thma}

\subsection{Jacobian algebra}

\begin{dfa}
 Let $Q$ be a finite quiver and $W$ a potential on $Q$. The \emph{Jacobian algebra} $J(Q,W)$ is the zeroth homology of the Ginzburg algebra $\Gamma(Q,W)$. This is the quotient algebra $$kQ/\langle \partial_a W, a\in Q_1\rangle$$
where $\langle \partial_a W, a\in Q_1\rangle$ is the two-sided ideal generated by the $\partial_aW$.
\end{dfa}
Remark: We follow the terminology of H.~Derksen, J.~Weyman and A.~Zelevinsky (\cite{DWZ} definition 3.1).

In recent works, B. Keller \cite{Kel10} and A. Buan, O. Iyama, I. Reiten and D. Smith \cite{BIRSm} have shown independently the following result:

\begin{thma}[Keller, Buan-Iyama-Reiten-Smith]
Let $T$ be a cluster-tilting object in the cluster category $\Cc_Q$ associated to an acyclic quiver $Q$.  Then there exists a quiver with potential $(Q',W)$ such that $\End_{\Cc_Q}(T)$ is isomorphic to $J(Q',W)$. 
\end{thma}

\subsection{Jacobi-finite quiver with potentials}

The quiver with potential $(Q,W)$ is called \emph{Jacobi-finite} if the Jacobian algebra $J(Q,W)$ is finite-dimensional.

\begin{dfa}
 Let $(Q,W)$ be a Jacobi-finite quiver with potential. Denote by $\Gamma$ the Ginzburg dg algebra $\Gamma(Q,W)$. Let $\per \Gamma$ be the thick subcategory of $\Dd\Gamma$ generated by $\Gamma$ and $\Dd^b\Gamma$ the full subcategory of $\Dd\Gamma$ of the dg $\Gamma$-modules whose homology is of finite total  dimension. The \emph{cluster category} $\Cc_{(Q,W)}$ associated to $(Q,W)$ is defined as the quotient of triangulated categories $\per\Gamma/\Dd^b\Gamma.$ 
\end{dfa}

Combining theorem \ref{clustertilting} and theorem \ref{ginzburg3CY} we get the result:

\begin{thma}\label{potentiel2CY}
 Let $(Q,W)$ be a Jacobi-finite quiver with potential. The cluster category $\Cc_{(Q,W)}$ associated to $(Q,W)$ is $\Hom$-finite and $2$-CY. Moreover the image $T$ of the free module $\Gamma$ in the quotient $\per\Gamma/\Dd^b\Gamma$ is a cluster-tilting object. Its endomorphim algebra is isomorphic to the Jacobian algebra $J(Q,W)$.
\end{thma}

As a direct consequence of this theorem we get the corollary:

\begin{cora}\label{CYtilted}
Each finite-dimensional Jacobi algebra $\Jj(Q,W)$ is
$2$-CY-tilted in the sense of I. Reiten (cf. \cite{Rei2}),
i.e. it is the endomorphism algebra of some cluster-tilting object of a
$2$-CY category.
\end{cora}

\begin{dfa}
 Let $(Q,W)$ and $(Q',W')$ be two quivers with potential. A \emph{triangular extension} between $(Q,W)$ and $(Q',W')$ is a quiver with potential $(\bar{Q},\bar{W})$ where
\begin{itemize}
 \item $\bar{Q}_0=Q_0\cup Q_0'$;
\item $\bar{Q}_1=Q_1\cup Q_1'\cup\{a_i,i\in I\}$, where for each $i$ in the finite index set $I$, the source of $a_i$ is in $Q_0$ and the tail of $a_i$ is in $Q'_0$;
\item $\bar{W}=W+W'$.
\end{itemize}

\end{dfa}

\begin{prop}\label{classedeJ} 
Denote by $\Jj\Ff$ the class of Jacobi-finite quivers with potential. The class $\Jj\Ff$ satisfies the properties:
 \begin{enumerate}
  \item it contains all acyclic quivers (with potential 0);
\item   it is stable under mutation of quivers with potential defined in \cite{DWZ};
\item it is stable under triangular extensions.
 \end{enumerate}

\end{prop}

\begin{proof}
 \begin{enumerate}
  \item This is obvious since the Jacobi algebra $J(Q,0)$ is isomorphic to $kQ$.
\item This is corollary 6.6 of \cite{DWZ}.
\item Let $(Q,W)$ and $(Q',W')$ be two quivers with potential in $\Jj\Ff$ and $(\bar{Q},\bar{W})$ a triangular extension. Let $\bar{Q}_1=Q_1\cup Q'_1\cup F$ be the set of arrows of $\bar{Q}$. We have then
$$k\bar{Q}= kQ'\ten_{R'}(R'\oplus kF\oplus R)\ten_RkQ$$ where $R$ is the semi-simple algebra $kQ_0$ and $R'$ is $kQ_0'$. Let $\bar{W}$ be the potential $W+W'$ associated to the triangular extension. If $a$ is in $Q_1$, then $\partial_a\bar{W}=\partial_aW$, if $a$ is in $Q_1'$ then  $\partial_a\bar{W}=\partial_aW'$ and if $a$ is in $F$, then $\partial_a\bar{W}=0$. Thus we have isomorphisms
\begin{eqnarray*}
 J(\bar{Q},\bar{W})& = &k\bar{Q}/\langle \partial_a\bar{W}, a\in \bar{Q}_1\rangle\\ & \simeq & kQ'\ten_{R'}(R'\oplus kF\oplus R)\ten_RkQ/\langle \partial_aW, a\in Q_1, \partial_bW', b\in Q'_1\rangle\\
&\simeq & kQ'/\langle \partial_bW', b\in Q'_1\rangle\ten_{R'}(R'\oplus kF \oplus R)\ten_RkQ/\langle \partial_aW, a\in Q_1\rangle\\
&\simeq & J(Q',W')\ten_{R'}(R'\oplus kF\oplus R)\ten_{R}J(Q,W).
\end{eqnarray*}
Thus if $J(Q',W')$ and $J(Q,W)$ are finite-dimensional, $J(\bar{Q},\bar{W})$ is finite-dimensional since $F$ is finite.
 \end{enumerate}

\end{proof}

In a recent work \cite{KY}, B. Keller and D. Yang proved the following:

\begin{thma}[Keller-Yang]
 Let $(Q,W)$ be a Jacobi-finite quiver with potential. Assume that $Q$ has no loops nor two-cycles. For each vertex $i$ of $Q$, there is a derived equivalence $$\Dd\Gamma(\mu_i(Q,W))\simeq \Dd\Gamma(Q,W),$$ where $\mu_i(Q,W)$ is the mutation of $(Q,W)$ at the vertex $i$ in the sense of \cite{DWZ}.
\end{thma}
Remark: in fact Keller and Yang proved this theorem in a more general setting. This also true if $(Q,W)$ is not Jacobi-finite, but then there is a derived equivalence between the completions of the Ginzburg dg algebras.

An other link between mutation of quivers with potential and mutations
of cluster-tilting objects is given in \cite{BIRSm} (theorem 5.1):

\begin{thma}[Buan-Iyama-Reiten-Smith]
Let $\Cc$ be a 2-CY triangulated category with a cluster-tilting
object $T$. If the endomorphism algebra $\End_\Cc(T)$ is isomorphic to
the Jacobian algebra $J(Q,W)$ for some quiver with potential $(Q,W)$,
and if no 2-cycles start in the vertex $i$ of $Q$, then we have an isomorphism
$$\End_\Cc(\mu_i(T))\simeq J(\mu_i(Q,W)).$$
\end{thma}

Combining these two theorems with theorem \ref{potentiel2CY}, we get the corollary:

\begin{cora}\label{corollaire}
 \begin{enumerate}
  \item If $Q$ is an acyclic quiver, and $W=0$, the cluster category $\Cc_{(Q,W)}$ is canonically equivalent to the cluster category $\Cc_Q$.
\item Let $Q$ be an acyclic quiver and $T$ a cluster-tilting object of $\Cc_Q$. If $(Q',W)$ is the quiver with potential associated  with the cluster-tilted algebra $\End_{\Cc_Q}(T)$ (cf. \cite{Kel10} \cite{BIRSm}),  then the cluster category $\Cc_{(Q,W)}$ is triangle equivalent to the cluster category $\Cc_{Q'}.$ 
 \end{enumerate}
\end{cora}

\begin{proof}
 \begin{enumerate}
  \item The cluster category $\Cc_{(Q,0)}$ is a $2$-CY category with a cluster-tilting object whose endomorphism algebra is isomorphic to $kQ$. Thus by \cite{Kel5}, this category is triangle equivalent to $\Cc_Q$.
\item In a cluster category, all cluster-tilting objects are mutation
  equivalent. Thus there exists a sequence of mutations which links
  $kQ$ to $T$. Moreover the quiver of a cluster-tilted algebra has no
  loops nor 2-cycles. Thus by theorem 5.1 of \cite{BIRSm}, the quiver
  with potential $(Q,W)$ is mutation equivalent to $(Q',0)$. Then the theorem of Keller and Yang \cite{KY} applies and we have an equivalence $$\Dd\Gamma(Q,W)\simeq \Dd\Gamma(Q',0).$$  Thus the categories $\Cc_{(Q,W)}$ and $\Cc_{(Q',0)}$ are triangle equivalent. By (1) we get the result. 
 \end{enumerate}

\end{proof}

\section{Cluster categories for non hereditary algebras}

\subsection{Definition and results of Keller}

Let $A$ be a finite-dimensional $k$-algebra of finite global dimension. The category $\Dd^bA$ admits a Serre functor $\nu_A=?\lten_A DA$ where $D$ is the duality $\Hom_k(?,k)$ over the ground field. The orbit category $$\Dd^bA/\nu_A\circ[-2]$$ is defined as follows:
\begin{itemize}
 \item the objects are the same as those of $\Dd^bA$;
\item if $X$ and $Y$ are in $\Dd^bA$ the space of morphisms is isomorphic  to the space $$\bigoplus_{i\in\mathbb{Z}}\Hom_{\Dd A}(X,(\nu_A^i Y[-2i]).$$
\end{itemize}
By Theorem 1 of \cite{Kel}, this category is triangulated if $A$ is derived equivalent to an hereditary category. This happens if $A$ is the endomorphism algebra of a tilting module over an hereditary algebra, or if $A$ is a canonical algebra (cf. \cite{Hap5}, \cite{Hap6}).

In general it is not triangulated and we define its \emph{triangulated hull} as the algebraic triangulated category $\Cc_A$ with the following universal property:
\begin{itemize}
 \item There exists an algebraic triangulated functor $\pi:\Dd^bA\rightarrow \Cc_A$.
\item Let $\Bb$ be a dg category and $X$ an object of $\Dd(A^{op}\ten \Bb)$. If there exists an isomorphism in $\Dd(A^{op}\ten \Bb)$ between $DA\lten_AX[-2]$ and $X$, then the triangulated algebraic functor $?\lten_AX:\Dd^bA\rightarrow \Dd\Bb$ factorizes through $\pi$.
\end{itemize}

Let $B$ be the dg algebra $A\oplus DA[-3]$. Denote by $p:B\rightarrow A$ the canonical projection. It induces a triangulated functor $p_*:\Dd^bA\rightarrow \Dd^bB$. Let $\langle A\rangle_B$ be the thick subcategory of $\Dd^bB$ generated by the image of $p_*$. By Theorem 2 of \cite{Kel} (cf. also \cite{Kelcorrection}), the triangulated hull of the orbit category $\Dd^bA/\nu_A\circ[-2]$ is the category $$\Cc_A=\langle A\rangle_B/\per B.$$ We call it the \emph{cluster category of $A$}. Note that if $A$ is the path algebra of an acyclic quiver, there is an equivalence $$\Cc_Q=\Dd^b(kQ)/\nu\circ[-2]\simeq \langle kQ\rangle_B/\per B.$$

\subsection{2-Calabi-Yau property}

The dg $B$-bimodule $DB$ is clearly isomorphic to $B[3]$, so it is not hard to check the following lemma:

\begin{lema}
 For each $X$ in $\per B$ and $Y$ in $\Dd^bB$ there is a functorial isomorphism
$$D\Hom_{\Dd B}(X,Y)\simeq \Hom_{\Dd B}(Y,X[3]).$$
\end{lema}

So we can apply results of section 1 and construct a bilinear bifunctorial form:

$$\beta_{XY}':\Hom_{\Cc_A}(X,Y)\times \Hom_{\Cc_A}(Y,X[2])\rightarrow k.$$

\begin{thma}\label{catama}
Let $X$ and $Y$ be objects in $\Dd=\Dd^bB$. If the spaces $\Hom_\Dd(X,Y)$ and $\Hom_\Dd(Y,X[3])$ are finite-dimensional, then the bilinear form
 $$\beta'_{XY}:\Hom_{\Cc_A}(X,Y)\times \Hom_{\Cc_A}(Y,X[2])\rightarrow
 k$$ is non-degenerate.
\end{thma}

Before proving this theorem, we recall some results about inverse limits of sequences of vector spaces that we will use in the proof. Let
$\xymatrix@-1.1pc{\ldots\ar[r] & V_p\ar[r]^\varphi & V_{p-1}\ar[r]^\varphi
  &\cdots\ar[r] & V_1\ar[r]^\varphi & V_0}$ be an inverse system of vector spaces (or vector space complexes) inverse system. We then have the following exact sequence

$$\xymatrix{0\ar[r] & V_\infty=\limproj V_p \ar[r] & \prod_p
  V_p \ar[r]^\Phi & \prod_q V_q \ar[r] & \limproj^1
V_p\ar[r] &0}$$
where $\Phi$ is defined by $\Phi(v_p)=v_p-\varphi(v_p)\in V_p\oplus
 V_{p-1}$ where $v_p$ is in $V_p$.

Recall two classical lemmas due to Mittag-Leffler:

\begin{lema}\label{l1}
If, for all $p$, the sequence of vector spaces  $W_i=\Imm(V_{p+i}\rightarrow V_p)$ is stationary, then $\limproj^1 V_p$ vanishes.
\end{lema}
This happens in particular when all vector spaces $V_p$ are finite-dimensional.

\begin{lema}\label{l2}
Let $\xymatrix{\ldots\ar[r] & V_p\ar[r]^\varphi & V_{p-1}\ar[r]^\varphi
  &\cdots\ar[r] & V_1\ar[r]^\varphi & V_0}$ be an inverse system of finite-dimensional vector spaces such that $V_\infty=\limproj V_p$ is also finite-dimensional. Let $V'_p$ be the image  of $V_\infty$ in $V_p$. The sequence $V'_p$ is stationary and we have $V'_\infty=\limproj V'_p= V_\infty$.
\end{lema}

\begin{proof}\emph{(of theorem \ref{catama})}
Let $X$ and $Y$ be objects of $\Dd^bB$ such that $\Hom_{\Dd^bB}(X,Y)$ is finite-dimensional. We will prove that there exists a local $\per B$-cover of $X$ relative to $Y$.

Let $\xymatrix{P_\bullet:\ldots\ar[r] & P_{n+1}\ar[r] & P_n\ar[r] & P_{n-1}\ar[r] & \ldots\ar[r]
   & P_0}$ be a projective resolution of $X$. The complex $P_\bullet$ has components in $\per B$, and its homology vanishes in all degrees except in degree zero, where it is $X$.
 Let $P_{\leq n}$ and $P_{>n}$ be the natural truncations, and denote by $Tot(P)$ the total complex associated to $P_\bullet$.
 For all $n\in\mathbb{N}$, there is an exact sequence of dg $B$-modules:
$$\xymatrix{0\ar[r] & Tot(P_{\leq n})\ar[r] & Tot(P)\ar[r] & Tot(P_{> n})\ar[r] & 0}$$
The complex $Tot(P)$ is quasi-isomorphic to $X$, and the complex $Tot(P_{\leq n})$ is in $\per B$. Moreover, $Tot (P)$ is the colimit of $ Tot(P_{\leq n})$. Thus by definition, we have the following equalities
$$\begin{array}{rcl} \Hh om ^\bullet _B(Tot (P),Y) & = & \Hh om ^\bullet
  _B( \underset{\rightarrow}{\rm{colim}} Tot (P_{\leq n}), Y)\\ & = &
  \limproj  \Hh om ^\bullet _B(Tot (P_{\leq n}),Y).\end{array}$$

Denote by $V_p$ the complex $\Hh om ^\bullet _B(Tot (P_{\leq p}),Y)$. In the inverse system  $$\xymatrix{\ldots\ar[r] & V_p\ar[r]^\varphi & V_{p-1}\ar[r]^\varphi
  &\cdots\ar[r] & V_1\ar[r]^\varphi & V_0},$$ all the maps are surjective, so by lemma \ref{l1}, there is a short exact sequence
$$\xymatrix{0\ar[r] & V_\infty\ar[r] & \prod_p
  V_p \ar[r]^\Phi & \prod_q V_q \ar[r] & 0}$$ which induces a long exact sequence in cohomology
$$\xymatrix@-1.1pc{\cdots\prod_q H^{-1}V_q\ar@{->>}[dr] \ar[rr] & &
  H^0(V_\infty)\ar@{->>}[dr] \ar[rr] & & \prod H^0V_p \ar[r] &
  \prod H^0V_q\cdots  \\   & \limproj^1 H^{-1}V_p
  \ar@{>->}[ur] & &  \limproj H^{0}V_p
  \ar@{>->}[ur] &  &}.$$
We have the equalities
 $$\begin{array}{rcl} H^0(V_\infty) & = & H^0(\Hh om ^\bullet
  _B(Tot (P),Y))\\ & = &\Hom_\Hh(Tot(P),Y)\\ & = &\Hom_\Dd(X,Y).\end{array}$$

Denote by $W_p$ the complex $\Hom_\Dd(Tot(P_{\leq p}),Y)$ and by $U_p$ the complex $H^{-1}(V_p)=\Hom_\Dd(Tot(P_{\leq p}),Y[-1])$.
The spaces $(U_p)_p$ are finite-dimensional, so by lemma \ref{l1}, $\limproj^1 U_p$ vanishes and we have an isomorphism

$$H^0(\limproj V_p)=H^0(V_\infty)\simeq \limproj H^0(V_p).$$
The system $(W_p)_p$ satisfies the hypothesis of lemma \ref{l2}. In fact, for each integer $p$, the space $\Hom_\Dd(Tot(P_{\leq p}),Y)$ is finite-dimensional because $Tot(P_{\leq p})$ is in $\per B$. Moreover, by the last two equalities $W_\infty=\limproj W_p$ is isomorphic to $\Hom_\Dd(X,Y)$ which is finite-dimensional by hypothesis. By lemma \ref{l2}, the system $(W'_p)_p$ formed by the image of $W_\infty$ in $W_p$ is stationary. More precisely, there exists an integer $n$ such that $W'_n=\limproj W'_p$. Moreover $W'_n$ is a subspace of $W_n=\Hom_\Dd(Tot(P_{\leq n}),Y)$ and there is an injection $$\xymatrix{\Hom_\Dd(X,Y)\ar@{^(->}[r] &
  \Hom_\Dd(Tot(P_{\leq n}),Y)}.$$ This yields a local $\per B$-cover of $X$ relative to $Y$.

The spaces $\Hom_\Dd(N,X)$ and $\Hom_\Dd(X,N)$ are finite-dimensional for $N$ in $\per B$ and $X$ in $\Dd^bB$. Thus by proposition \ref{restrictioncouverture}, there exists local $\per B$-envelopes. Therefore theorem \ref{nondegenere} applies and $\beta'$ is non-degenerate.
 
\end{proof}
 
\begin{cora}\label{2CY}
Let $A$ be a finite-dimensional $k$-algebra with finite global dimension. If the cluster category $\Cc_A$ is $\Hom$-finite, then it is $2$-CY as a triangulated category.
\end{cora}

\begin{proof}
 Denote by $p_*:\Dd^bA\rightarrow \Dd^bB$ the restriction of the projection $p:B\rightarrow A$.

Let $X$ and $Y$ be in $\Dd^b(A)$.  By hypothesis, the vector spaces 

$$\bigoplus_{p\in\mathbb{Z}}\Hom_{\Dd^bA}(X, \nu_A^p Y[-2p]) \quad
\text{and} \quad \bigoplus_{p\in\mathbb{Z}}\Hom_{\Dd^bA}(Y, \nu_A^p
X[-2p+3])$$ are finite-dimensional. But by \cite{Kel}, the space $\Hom_{\Dd^bB}(p_*X,p_*Y)$ is isomorphic to $$\bigoplus_{p\geq 0}\Hom_{\Dd^bA}(X, \nu_A^p Y[-2p]),$$ so is finite-dimensional. For the same reasons, the space $\Hom_{\Dd^bB}(Y,X[3])$ is also finite-dimensional. Applying theorem \ref{catama}, we get a non-degenerate bilinear form $\beta'_{p_*X,p_*Y}$.
The non-degeneracy property is extension closed, so for each $M$ and $N$ in $\langle A\rangle_B$, the form $\beta'_{MN}$ is non-degenerate.

\end{proof}

\subsection{Case of global dimension 2}\label{dimensionglobale2}

In this section we assume that $A$ is a finite-dimensional $k$-algebra of global dimension $\leq 2$. 

\subsubsection*{Criterion for $\Hom$-finiteness}

The canonical $t$-structure on the derived category $\Dd=\Dd^bA$ satisfies the property: 
\begin{lema}\label{lemmenu}
We have the following inclusions $\nu(\Dd_{\geq 0})\subset \Dd_{\geq -2}$ and $\nu^{-1}(\Dd_{\leq 0})\subset \Dd_{\leq 2}$. Moreover, the space $\Hom_\Dd(U,V)$ vanishes for all $U$ in $\Dd_{\geq 0}$ and all $V$ in $\Dd_{\leq -3}$.
\end{lema}

\begin{prop}\label{endomorphisme}
Let $X$ be the $A$-$A$-bimodule $\Ext^2_A(DA,A)$. The endomorphism algebra $\widetilde{A}=\End_{\Cc_A}(A)$ is isomorphic to the tensor algebra $T_A X$ of $X$ over $A$.
\end{prop}

\begin{proof}
By definition, the endomorphism space $\End_{\Cc_A}(A)$ is isomorphic to $$\bigoplus_{p\in\mathbb{Z}}\Hom_\Dd(A,\nu^p
A[-2p])$$
For $p\geq1$, the object $\nu^pA[-2p]$ is in $\Dd_{\geq 2}$ since $\nu A$ is in $\Dd_{\geq 0}$. So since $A$ is in $\Dd_{\leq 0}$, the space $\Hom_\Dd(A,\nu^p
A[-2p])$ vanishes.

The functor $\nu=?\lten_ADA$ admits an inverse $\nu^{-1}=-\lten_A ~R\Hh om_A(DA,A).$
Since the global dimension of $A$ is $\leq 2$, the homology of the complex $R\Hh om_A(DA,A)$ is concentrated in degrees $0$, $1$ and $2$ :
\begin{eqnarray*}
 H^0(R\Hh om_A(DA,A))&=&\Hom_\Dd(DA,A)\\
H^1(R\Hh om_A(DA,A))&=&\Ext^1_A(DA,A)\\
H^2(R\Hh om_A(DA,A))&=&\Ext^2_A(DA,A).
\end{eqnarray*}

Let us denote by $Y$ the complex  $R\Hh om_A(DA,A)[2]$.  We then have 
$$\nu^{-p}A[2p]=A\lten_A (Y^{\lten_A p})=Y^{\lten_A p}.$$ Therefore we get the following equalities
\begin{eqnarray*}\Hom_{\Dd A}(A,S^{-p} A[-2p])& = & \Hom_{\Dd
    A}(A,Y^{\lten_A p})\\
 & = & H^0(Y^{\lten_A p}).
\end{eqnarray*}
Since $H^0(Y)=X$, we conclude using the following easy lemma.
\end{proof}   
\begin{lema}\label{H0tensor}
Let $M$ and $N$ be two complexes of $A$-modules whose homology is concentrated in negative degrees. Then there is an isomorphism
$$H^0(M\lten_A N)\simeq H^0(M)\ten_A H^0(N).$$
\end{lema}

\begin{prop}\label{homfini}
 Let $A$ be a finite-dimensional algebra of global dimension $2$. The following properties are equivalent:
\begin{enumerate}
 \item the cluster category $\Cc_A$ is $\Hom$-finite;
\item the functor $?\ten_A\Ext^2(DA,A)$ is nilpotent;
\item the functor $\Tor_2^A(?,DA)$ is nilpotent;
\item there exists an integer $N$ such that there is an inclusion $\Phi^N(\Dd_{\geq 0})\subset \Dd_{\geq 1}$ where $\Phi$ is the autoequivalence $ \nu_A[-2]$ of the category $\Dd=\Dd^bA$ and $\Dd_{\geq 0}$ is the right aisle of the natural $t$-structure of $\Dd^bA$.
 \end{enumerate}

\end{prop}

\begin{proof}
$ 1\Rightarrow 2$: It is obvious by proposition \ref{endomorphisme}.

$2\Leftrightarrow 3\Leftrightarrow 4$: Let $\Phi$ be the autoequivalence $\nu_A[-2] $ of $\Dd^bA$. The functor $\Tor^2_A(?,DA)$ is isomorphic to $H^0\circ \Phi$ and $?\ten_A\Ext^2_A(DA,A)$ is isomorphic  to $H^0\circ \Phi^{-1}$. Thus it is sufficient to check that there are isomorphisms $$H^0\Phi\circ H^0\Phi\simeq H^0\Phi^2 \textrm{ and } H^0\Phi^{-1}\circ H^0\Phi^{-1}\simeq H^0\Phi^{-2}.$$
This is easy using Lemma \ref{H0tensor} since the algebra $A$ has global dimension $\leq 2$ .

$4\Rightarrow 1$: Suppose that there exists some $N\geq 0$ such that $\Phi^N(\Dd_{\geq 0})$ is included in $\Dd_{\leq 1}$.   For each object $X$ in $\Cc_A$, the class of the objects $Y$ such that the space $\Hom_{\Cc_A}(X,Y)$ (resp. $\Hom_{\Cc_A}(Y,X)$) is finite-dimensional, is extension closed. Therefore, it is sufficient to show that for all simples $S$, $S'$, and each integer $n$, the space $\Hom_{\Cc_A}(S,S'[n])$ is finite-dimensional.

There exists an integer  $p_0$ such that for all $p\geq p_0$ $\Phi^p(S')$ is in $\Dd_{\geq n
+1}$. Therefore, because of the defining properties of the $t$-structure, the space 
$$\bigoplus_{p\geq p_0}\Hom_\Dd(S,\Phi^p(S')[n]) $$ vanishes.
Similary, there exists an integer $q_0$ such that for all $q\geq q_0$, we have $\Phi^q(S)\in \Dd_{\geq -n
+3}$. Since the algebra $A$ is of global dimension $\leq 2$, the space $$\bigoplus_{q\geq q_0}
  \Hom_\Dd(\Phi^q(S),S'[n])$$ vanishes. Thus the space 
$$ \bigoplus_{p\in\mathbb{Z}}
  \Hom_\Dd(S,\Phi^p(S')[n])= \bigoplus_{p= -q_0}^{p_0}
  \Hom_\Dd(S,\Phi^p(S')[n])$$ is finite-dimensional.
\end{proof}

\subsubsection*{Cluster-tilting object}

In this section we prove the following theorem:
\begin{thma}\label{Aclustertilting}
 Let $A$ be a finite-dimensional $k$-algebra of global dimension $\leq 2$. If the functor $\Tor_2^A(?,DA)$ is nilpotent, then the cluster category $\Cc_A$ is $\Hom$-finite, $2$-CY and the object $A$ is a cluster-tilting object.
\end{thma}

We denote by $\Theta$ a cofibrant resolution of the dg $A$-bimodule $R\Hh om^\bullet_A(DA,A)$. Following \cite{Kel7} and \cite{Kel10}, we define the 3-derived preprojective algebra as the tensor algebra
$$\Pi_3(A)=T_A(\Theta[2]).$$ 
The complex $R\Hh om^\bullet_A(DA,A)[2]$ has its homology concentrated in degrees $-2$, $-1$ and $0$, and we have $$H^{-2}(\Theta[2])\simeq \Hom_{\Dd A}(DA,A),\ H^{-1}(\Theta[2])\simeq\Ext^1_A(DA,A)$$
$$\textrm{and}\ H^0(\Theta[2])\simeq \Ext^2_A(DA,A).$$ Thus the homology of the dg algebra $\Pi_3(A)$ vanishes in strictly positive degrees and we have
$$H^0\Pi_3A=T_A\Ext^2_A(DA,A)=\tilde{A}.$$
By proposition \ref{homfini}, it is finite-dimensional.
Moreover, Keller showed that $\Pi_3(A)$ is homologically smooth and bimodule 3-CY \cite{Kel10}. Thus we can apply theorem \ref{clustertilting} and we have the following result:

\begin{cora}
 The category $\Cc=\per \Pi_3A/\Dd^b\Pi_3A$ is $2$-CY and the free dg module $\Pi_3A$ is a cluster-tilting object in $\Cc$.
\end{cora}

To complete the proof of Theorem \ref{Aclustertilting} we now construct a triangle equivalence between $\Cc_A$ and $\Cc$ sending $A$ to $\Pi_3A$.

Let us recall a theorem of Keller (\cite{Kel9}, or theorem 8.5, p.96 \cite{tilting}):

\begin{thma}\label{kel10}[Keller]
 Let $B$ be a dg algebra, and $T$ an object of $\Dd B$. Denote by $C$ the dg algebra $R\Hh om^\bullet_B(T,T)$.
Denote by $\langle T \rangle_{B} $ the thick subcategory of $\Dd B$ generated by $T$. The functor $R\Hh om_B^\bullet(T,?):\Dd B\rightarrow \Dd C$ induces an algebraic triangle equivalence $$\xymatrix{R\Hh om_B^\bullet(T,?):\langle T \rangle_{B}\ar[r]^(.7)\sim & \per C.}$$ 
\end{thma}

Let us denote by $\Hh o(dgalg)$ the homotopy category of dg algebras, i.e. the localization of the category of dg algebras at the class of quasi-isomorphisms.
\begin{lema}
In $\Hh o(dgalg)$, there is an isomorphism between $\Pi_3 A$ and $R\Hom_B(A_B,A_B)$.
\end{lema}

\begin{proof}
 The dg algebra $B$ is $A\oplus (DA)[-3]$. Denote by $X$ a cofibrant resolution of the dg $A$-bimodule $DA[-2]$. Now look at the dg submodule of the bar resolution of $B$ seen as a bimodule over itself (see the proof of theorem 7.1 in \cite{Kel}):
$$\xymatrix{bar(X,B):\cdots \ar[r] & B\otimes_AX^{\otimes_A2}\otimes_AB\ar[r] & B\otimes_AX\otimes_AB\ar[r] & B\otimes_AB\ar[r] & 0}$$
This is a cofibrant resolution of the dg $B$-bimodule $B$.
Thus $A\otimes_Bbar(X,B)$ is a cofibrant resolution of the dg $B$-module $A$.
Therefore, we have the following isomorphisms
\begin{eqnarray*}
 R\Hh om_B^\bullet(A_B,A_B) &\simeq & \Hh om^\bullet_B(A\otimes_Bbar(X,B),A)\\ & \simeq & \prod_{n\geq 0} \Hh om^\bullet_B(A\otimes_A X^{\otimes_A n}\otimes_A B,A_B) \\ & \simeq & \prod_{n\geq 0} \Hh om_A^\bullet (X^{\otimes_A n},\Hom_B(B,A_B)_A)\\ & \simeq & \prod_{n\geq 0} \Hh om_A^\bullet (X^{\otimes_A n},A_A),
\end{eqnarray*}
where the differential on the last complex is induced by that of $X^{\otimes_A n}$.
Note that \begin{eqnarray*}\Hh om^\bullet_A(X,A) &=&R\Hh om^\bullet_A(DA[-2],A)\\
               &=& R\Hh om^\bullet_A(DA,A)[2] = \Theta[2].
              \end{eqnarray*}
We can now use the following lemma:
\begin{lema}
 Let $A$ be a dg algebra, and $L$ and $M$ dg $A$-bimodules such that $M_A$ is perfect as right dg $A$-module. There is an isomorphism in $\Dd(A^{op}\ten A)$ $$R\Hh om^\bullet_A(L,A)\lten_AR\Hh om^\bullet_A(M,A) \simeq R\Hh om^\bullet_A(M\lten_AL,A).$$
\end{lema}

\begin{proof}
Let $X$ and $M$ be dg $A$-bimodules. The following morphism of $\Dd(A^{op}\ten A)$
\begin{eqnarray*}
 X\lten_AR\Hh om_A(M,A) & \longrightarrow & R\Hh om_A(M,X)\\
x\otimes \varphi & \mapsto & (m\mapsto x\varphi(m))
\end{eqnarray*}
is clearly an isomorphism for $M=A$. Thus it is an isomorphism if $M$ is perfect as a right dg $A$-module.
Applying this to  the right dg $A$-module $R\Hh om_A(L,A)$, we get an isomorphism of dg $A$-bimodules
$$R\Hh om_A(L,A)\lten_AR\Hh om_A(M,A) \simeq R\Hh om_A(M,R\Hh om_A(L,A)).$$
Finally, by adjunction we get an isomorphism of dg $A$-bimodules
$$R\Hh om_A(L,A)\lten_AR\Hh om_A(M,A) \simeq R\Hh om_A(M\lten_A L,A).$$
\end{proof}
Therefore, the dg $A$-bimodule $\Hh om_A^\bullet (X^{\otimes_A n},A_A)$ is isomorphic to $(\Theta[2])^{\otimes_A n}$, and there is an isomorphism of dg algebras $$R\Hh om_B^\bullet(A_B,A_B)\simeq \bigoplus_{n\geq 0}(\theta[2])^{\lten_A n}=\Pi_3(A)$$
because for each $p\in \mathbb{Z}$, the group $H^p(\theta[2]^{\lten_A n})$ vanishes for all $n\gg0$.
\end{proof}

By theorem \ref{kel10}, the functor $R\Hh om ^\bullet_B(A_B,?)$ induces an equivalence between the thick subcategory $\langle A\rangle_B$ of $\Dd B$ generated by $A$, and $\per \Pi_3(A)$. Thus  we get a triangle equivalence that we will denote by $F$:
$$\xymatrix{F=R\Hh om^\bullet_B(A_B,?):\langle A\rangle_B \ar[rr]^(.7)\sim && \per \Pi_3A}$$  
This functor sends the object $A_B$ of $\Dd^bB$ onto the free module $\Pi_3A$ and the free $B$-module $B$ onto $R\Hh om ^\bullet_B(A_B,B)\simeq R\Hh om^\bullet _B(A_B,DB[-3])$, that is to say onto $(DA)[-3]_{\Pi_3A}$. So $F$ induces an equivalence
$$\xymatrix{F:\per B=\langle B \rangle_B\ar[r]^(.4)\sim & \langle DA[-3]\rangle_{\Pi_3A} =\langle A\rangle_{\Pi_3A}.}$$ 

\begin{lema}
 The thick subcategory $\langle A\rangle_{\Pi_3A}$ of $\Dd\Pi_3A$ generated by $A$  is $\Dd^b\Pi_3 A$.
\end{lema}

\begin{proof}
The algebra $A$ is finite-dimensional, thus $\langle A\rangle_{\Pi_3A}$ is obviously included in $\Dd^b\Pi_3 A$. Moreover, the category $\Dd^b\Pi_3 A$ equals $\langle \mod H^0(\Pi_3A)\rangle_{\Pi_3A}$ by the existence of the $t$-structure. The dg algebra $\Pi_3 A$ is the tensor algebra $T_A(\theta[2])$ thus there is a canonical projection $\Pi_3A\rightarrow A$ which yields a restriction functor $\Dd^bA\rightarrow \Dd^b(\Pi_3A)$ respecting the $t$-structure:
$$\xymatrix{\mod H^0\Pi_3 A=\Hh\ar@{^(->}[r]& \Dd^b(\Pi_3A)\\ \mod A\ar[u]\ar@{^(->}[r] & \Dd^bA\ar[u]}$$
This restriction functor induces a bijection in the set of isomorphism classes of simple modules because the kernel of the map $H^0(\Pi_3A)\rightarrow A$ is a nilpotent ideal (namely the sum of the tensor powers over $A$ of the bimodule $\Ext_A^2(DA,A)$). Thus each simple of $\mod H^0\Pi_3A$ is in $\langle A\rangle_{\Pi_3A}$ and we have
$$\langle A\rangle_{\Pi_3 A}\simeq  \langle \mod H^0(\Pi_3A)\rangle_{\Pi_3A}\simeq \Dd^b\Pi_3 A.$$
 
\end{proof}

\begin{proof}\textit{(of theorem \ref{Aclustertilting})}
 By proposition \ref{homfini} and corollary \ref{2CY}, the cluster category is $\Hom$-finite and $2$-CY. 
Furthermore, the functor $F=R\Hh om^\bullet_B(A_B,?)$ induces the following commutative square:
$$\xymatrix{F:\langle A\rangle_B\ar[r]^\sim & \per \Pi_3A\\ \per B\ar[r]^\sim\ar@{^(->}[u] &\Dd^b\Pi_3A\ar@{^(->}[u]}$$
Thus $F$ induces a triangle equivalence
$$\xymatrix{\Cc_A=\langle A \rangle_B/\per B\ar[rr]^\sim && \per \Pi_3A/\Dd^b\Pi_3A=\Cc}$$
sending the object $A$ onto the free module $\Pi_3A$. By theorem \ref{clustertilting}, $A$ is therefore a cluster-tilting object of the cluster category $\Cc_A$.
\end{proof}

\subsubsection*{Quiver of the endomorphism algebra of the cluster-tilting object}

 Let $A=kQ/I$ be a finite-dimensional $k$-algebra of global dimension $\leq 2$. Suppose that $I$ is an admissible ideal generated by a finite set of minimal relations $r_i$, $i\in J$ where for each $i\in J$, the relation $r_i$ starts at the vertex $s(r_i)$ and ends at the vertex $t(r_i)$.  Let $\widetilde{Q}$ be the following quiver:
\begin{itemize}
 \item the set of the vertices of $\widetilde{Q}$ equals that of $Q$;
\item the set of arrows of $\widetilde{Q}$ is obtained from that of $Q$ by adding a new arrow $\rho_i$ with source $t(r_i)$ and target $s(r_i)$ for each $i$ in $J$.
\end{itemize}

We then have  the following proposition, which has essentially been proved by I. Assem, T. Br{\"u}stle and R. Schiffler \cite{Ass} (thm 2.6). The proposition is also proved in \cite{Kel10}.

\begin{prop}\label{carquoisdetildeA}
If the algebra $\End_{\Cc_A}(A)=\widetilde{A}$ is finite-dimensional, then
 its quiver is $\widetilde{Q}$.
\end{prop}

\begin{proof}
Let $B$ be a finite-dimensional algebra. The vertices of its quiver are determined by the quotient $B/rad(B)$ and the arrows are determined by $rad(B)/rad^2(B)$. Denote by $X$ the $A$-$A$-bimodule $\Ext^2_A(DA,A)$. Since $X\ten_AX$ is in $rad^2(B)$, the quiver of $\tilde{A}=T_AX$ is the same as the quiver of the algebra $A\rtimes X$.  The proof is then exactly the same as in \cite{Ass} (thm 2.6).
 
\end{proof}

\subsubsection*{Example} We refer to \cite{Gei} for this example.
Let $Q$ be a Dynkin quiver. Let $A$ be its Auslander algebra. The algebra $A$ is of global dimension $\leq 2$. The category $\mod A$ is equivalent to the category $\mod (\mod kQ)$ of finitely presented functors $(\mod kQ)^{op}\rightarrow \mod k$. The projective indecomposables of $\mod A$ are the representable functors $U^\wedge=\Hom_{kQ}(?,U)$ where $U$ is an indecomposable $kQ$-module. Let $S$ be a simple $A$-module. Since $A$ is finite-dimensional, this simple is associated to an indecomposable $U$ of $\mod kQ$. If $U$ is not projective, then it is easy to check that in $\Dd^b(A)$ the simple $S_U$ is isomorphic to the complex:
$$\xymatrix{\cdots\ar@<1ex>[r] & \underset{-3}{0}\ar@<1ex>[r] & \underset{-2}{(\tau U)^\wedge}\ar@<1ex>[r] &
 \underset{-1}{ E^\wedge}\ar@<1ex>[r] &  \underset{0}{U^\wedge}\ar@<1ex>[r] & \underset{1}{0}\ar@<1ex>[r] &\cdots}$$
where $\xymatrix{0\ar[r] & \tau U\ar[r] & E\ar[r] & U\ar[r] & 0}$ is the Auslander-Reiten sequence associated to $U$. Thus $\Phi (S_U)=\nu S_U[-2]$ is the complex:
$$\xymatrix{\cdots\ar@<1ex>[r] & \underset{-1}{0}\ar@<1ex>[r] & \underset{0}{(\tau U)^\vee}\ar@<1ex>[r] &
 \underset{1}{ E^\vee}\ar@<1ex>[r] &  \underset{2}{U^\vee}\ar@<1ex>[r] & \underset{3}{0}\ar@<1ex>[r] &\cdots}$$ where $U^\vee$ is the injective $A$-module $D\Hom_{kQ}(U,?)$.
It follows from the Auslander-Reiten formula that this complex is quasi-isomorphic to the simple $S_{\tau U}$.

If $U$ is projective, then $S_U$ is isomorphic in $\Dd^b(A)$ to  $$\xymatrix{\cdots \ar@<1ex>[r] & \underset{-2}{0}\ar@<1ex>[r] &
  \underset{-1}{(rad U)^\wedge}\ar@<1ex>[r] & \underset{0}{U^\wedge}\ar@<1ex>[r] &
  \underset{1}{0} \ar@<1ex>[r]& \cdots},$$ and then  $\Phi(S_U)$ is in $\Dd_{\geq 1}$.
Since for each indecomposable $U$ there is some $N$ such that $\tau^NU$ is projective, there is some $M$ such that $\Phi^M(\Dd_{\geq 0})$ is included in $\Dd_{\geq 1}$.  By proposition \ref{homfini}, the cluster category $\Cc_A$ is $\Hom$-finite, and $2$-CY by corollary \ref{2CY}.

The quiver of $A$ is the Auslander-Reiten quiver of $\mod kQ$. The minimal relations of the algebra $A$ are given by the mesh relations. Thus the quiver of $\tilde{A}$ is the same as that of $A$ in which arrows $\tau x\rightarrow x$ are added for each non projective indecomposable $x$.

For instance, if $Q$ is $A_4$ with the orientation $\xymatrix{1\ar[r] & 2\ar[r] & 3\ar[r] & 4}$, then the quiver of the algebra $\tilde{A}$ is the following $$
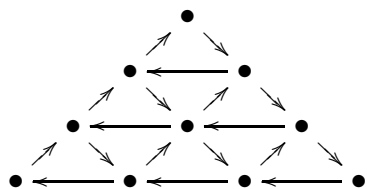
\xymatrix@-1.3pc{&&&\bullet\ar[dr]&&&\\&&\bullet\ar[ur]\ar[dr]&&\bullet\ar[ll]\ar[dr]&&
  \\&\bullet\ar[ur]\ar[dr]&&\bullet\ar[ll]\ar[ur]\ar[dr]&&\bullet\ar[ll]\ar[dr]&\\ \bullet\ar[ur]&&\bullet\ar[ll]\ar[ur]&&\bullet\ar[ll]\ar[ur]&&\bullet\ar[ll].}$$

\section{Stable module categories as cluster categories}

\subsection{Definition and first properties }

Let $B$ be a \emph{concealed} algebra \cite{Rin}, \emph{i.e.} the endomorphism algebra of a preinjective tilting module over a finite-dimensional hereditary algebra. Let $H$ be a postprojective slice of $\mod B$. We denote by $add(H)$ the smallest subcategory of $\mod B$ which contains $H$ and which is stable under taking direct summands.
 Let $Q$ be the quiver such that $\End_{B}(H)$ is the path algebra $kQ$ and let $Q_0=\{1,\cdots,n\}$ be its set of vertices.
By Happel \cite{Hap2}, there is a triangle equivalence:
   $$\xymatrix{\Dd^b(B)\ar@<.5ex>[rrr]^{DR\Hom_{B}(?,H)} &&& \Dd^b(kQ)\ar@<.5ex>[lll]^{(D?)\lten_{kQ} H}.}$$

Notice that these functors induce quasi-inverse equivalences between $add(H)$ and the subcategory of finite-dimensional injective $kQ$-modules.

Define $\Mm$ as the following subcategory of $\mod B$:
$$ \Mm
  =\{X\in \mod B \ | \ \Ext^1_{B}(X,H)=0\}=\{ X\in\mod B\ | \  X \textrm{ is cogenerated by }H\}$$ 
We denote by $\tau_B$ the AR-translation of the category $\mod B$ and by $\tau_{\Dd}$ the AR-translation of $\Dd^b B$. 

The following proposition is a classical result in tilting theory (see for example \cite{Rin}).

\begin{prop}\label{propbase}
\begin{enumerate}
 \item For each $X$ in $\Mm$ there exists a triangle $$\xymatrix{X\ar[r] & H_0\ar[r] & H_1\ar[r]
 & X[1]} $$ in $\Dd^b(\mod B)$
 functorial in $X$  with $H_0$ and $H_1$ in $add(H)$;
\item $\Mm\subset\mod B$ is closed under kernels so in particular, $\Mm$ is closed under $\tau_B$;
\item
for each indecomposable $X$ in $\Mm$ there exists a unique $q\geq 0$ such that $\tau^{-q}_B X$
 is in $add(H)$;
 \item the category $\Mm$ has finitely
 many indecomposables.
\end{enumerate}

\end{prop}

\subsubsection*{Hom-finiteness}

Let $\Mmm$ be the quotient $\Mm/add(H)$. Denote by $p:\Mm\rightarrow \Mmm$ the canonical projection. Since $H$ is a slice, we have the following properties.
\begin{prop}\label{propmbar}
\begin{enumerate}
 \item The category $\Mmm$ is equivalent to the full subcategory of $\Mm$ whose objects do not have non zero direct factors in $add(H)$. We denote by $i:\Mmm\rightarrow \Mm$ the associated inclusion.
\item The category $\Mmm\subset\mod B$ is closed under kernels, and hence under $\tau_B$.
\item The right exact functor $i:\mod\Mmm\rightarrow\mod\Mm$ induced by $i:\Mmm\rightarrow \Mm$ is isomorphic to the restriction along $p$.
\end{enumerate}
\end{prop}

\begin{prop}
Let $A$ be the endomorphism algebra
$\End_B(\bigoplus_{M\in\ind\Mmm}M)$. The global dimension of $A$
is at most $2$.
\end{prop}

\begin{proof}
 There is an equivalence of categories between $\mod A$ and $\mod \Mmm$. Since $\Mmm$ is stable under kernels, the global dimension of $A$ is $\leq 2$.
\end{proof}

\begin{thma}\label{homfini2}
The cluster category $\Cc_A$ is a $\Hom$-finite, $2$-CY category, and the object $A$ is a cluster-tilting object in $\Cc_A$.
\end{thma}

\begin{proof}
 Using corollary \ref{2CY} and  theorem \ref{homfini}, we just have to check that the functor $\Tor^2_A(?,DA)$ is nilpotent. Since there are finitely many indecomposables in $\Mmm$, the proof is the same as for an Auslander algebra (cf. the examples of section \ref{dimensionglobale2}).
\end{proof}

\subsubsection*{Construction of the functor $F:\mod \Mm \rightarrow
  \textrm{f.l.}\Lambda$}\label{construction}

Denote by $\Ii(kQ)$ the subcategory of the preinjective modules of $\mod kQ$.

\begin{prop}\label{constructionP}
 There exists a $k$-linear functor $P:\Ii(kQ)\rightarrow \Mm$ unique up to isomorphism such that
\begin{itemize}
 \item $P$ restricted to subcategory of the injective $kQ$-modules is isomorphic to the restriction of the functor $D(?)\ten_{kQ}H$;
\item for each indecomposable $X$ in $\Ii(kQ)$ such that $P(X)$ is not projective, the image $$\xymatrix{0\ar[r] & P(\tau_\Dd X)\ar[r]^{Pi} & P(E)\ar[r]^{Pp}&P(X)\ar[r] & 0}$$ of an Auslander-Reiten sequence in $\mod kQ$ ending at $X$ $$\xymatrix{0\ar[r] & \tau_\Dd X\ar[r]^i & E\ar[r]^p&X\ar[r] & 0}$$ is an Auslander-Reiten sequence in $\mod B$ ending at $P(X)$.
\end{itemize}
 Moreover, the functor $P$ is full, essentially surjective, and satisfies $P\circ\tau_\Dd\simeq \tau_B\circ P$.
\end{prop}

\begin{proof}
 The Auslander-Reiten quivers $\Gamma_\Ii$ of $\Ii(kQ)$ and $\Gamma_\Mm$ of $\Mm$ are connected translation quivers. Each vertex of $\Gamma_\Ii$ is of the form $\tau_\Dd^qx$ with $q\geq 0$ and $x$ indecomposable injective. Each vertex of $\Gamma_\Mm$ is of the form $\tau_B^qx$ where $x$ is in $add(H)$ ((3) of proposition \ref{propbase}). Moreover, there is a canonical isomorphism of quivers $\bar{P}: \Gamma_{DkQ}\rightarrow \Gamma_{add(H)}$. Thus we can inductively construct a morphism of quivers (that we will still denote by $\bar{P}$) $\bar{P}:\Gamma_\Ii\rightarrow \Gamma_\Mm$ extending $\bar{P}$ such that:
\begin{itemize}
 \item $\bar{P}(\tau_\Dd x)=\tau_B\bar{P}(x)$ for each vertex $x$ of $\Gamma_\Ii$ ;
\item $\bar{P}(\sigma_\Dd\alpha)=\sigma_B\bar{P}(\alpha)$ for each arrow $\alpha:x\rightarrow y$ of $\Gamma_\Ii$, where $\sigma_\Dd \alpha$ (resp. $\sigma_B \beta$) denotes the arrow $\tau_\Dd y\rightarrow x$ (resp. $\tau_B y\rightarrow x$) such that the mesh relations in $\Gamma_\Ii$ (resp. in $\Gamma_\Mm$) are of the form $\sum_{t(\alpha)=x}\sigma_\Dd(\alpha)\alpha$ (resp. $\sum_{t(\beta)=x}\sigma_B(\beta)\beta$).
\end{itemize}
Clearly, this morphism of translation quivers induces surjections in the sets of vertices and the sets of arrows.

The categories $\Ii(kQ)$ and $\Mm$ are standard, \emph{i.e.} $k$-linearly equivalent to the mesh categories of their Auslander-Reiten quivers. Up to isomorphism, an equivalence $k(\Gamma_\Ii)\rightarrow \Ii(kQ)$ is uniquely determined by its restriction to a slice. Thus there exists a $k$-linear functor $P:\Ii(kQ)\rightarrow \Mm$ unique up to isomorphism which is equal to $D(?)\ten_{kQ}H$ on the slice of the injectives and such that the square 
$$\xymatrix{k(\Gamma_{\Ii})\ar[d]^{\bar{P}}\ar[r]^\sim & \Ii(kQ)\ar[d]^P\\k(\Gamma_\Mm)\ar[r]^\sim & \Mm}$$
is commutative. This functor $P$ sends Auslander-Reiten sequences $$\xymatrix{0\ar[r] & \tau_\Dd X\ar[r]^i & E\ar[r]^p&X\ar[r] & 0}$$ to Auslander-Reiten sequences $$\xymatrix{0\ar[r] & \tau_B P(X)\ar[r]^{Pi} & P(E)\ar[r]^{Pp}&P(X)\ar[r] & 0}$$ if $P(X)$ is not projective. Since $\bar{P}$ is surjective, $P$ is full and essentially surjective.
\end{proof}

\begin{lema}\label{kernel}
 Let $X$ and $Y$ be indecomposables in $\Ii(kQ)$. The kernel of the map $\Hom_{kQ}(X,Y)\rightarrow \Hom_B(PX,PY)$ is generated by compositions of the form $X\rightarrow Z\rightarrow Y$ where $Z$ is indecomposable and $P(Z)$ is zero.
\end{lema}
 \begin{proof}
If $P(X)$ or $P(Y)$ is zero this is obviously true. Suppose they are not. The mesh relations are minimal relations of the $k$-linear category $\Mm$ and $P$ is full. Thus the kernel of the functor $P$ is the ideal generated by the morphisms of the form 
$\xymatrix{U\ar[r]^g&V\ar[r]^h & W}$
where $\xymatrix{0\ar[r] & P(U)\ar[r]^{Pg} & P(V)\ar[r]^{Ph}&P(W)\ar[r] & 0}$ is an Auslander-Reiten sequence in $\Mm$. Since $P(U)$ is isomorphic to $\tau_B P(W)$, the indecomposable $U$ is isomorphic to $\tau_\Dd(W)$.  By the construction of $P$, $V$ is a direct factor of the middle term of the Auslander-Reiten sequence ending at $W$, and we can `complete'  the composition $\xymatrix{\tau_\Dd W\ar[r]^g & V\ar[r]^h &W}$ into an Auslander-Reiten sequence  $$\xymatrix{0\ar[r] & \tau_\Dd W\ar[rr]^-{\left({\bsm g\\g'\esm}\right)} & &V\oplus V'\ar[rr]^-{\left({\bsm h&h'\esm}\right)}&&W\ar[r] & 0}$$ with $P(V')=0$ and $P(g')=P(h')=0$. Thus the morpism $hg=-h'g'$  factors through an object in the kernel of $P$. 
 \end{proof}

 Now let $\Lambda$ be the preprojective algebra associated to the acyclic quiver $Q$. It is defined as the quotient $k\bar{Q}/(c)$ where $\bar{Q}$ is the double quiver of $Q$ which is obtained from $Q$ by adding to each arrow $a:i\rightarrow j$ an arrow $a^*:j\rightarrow i$ pointing in the opposite direction, and where $(c)$ is the ideal generated by the element $$c=\sum_{a\in Q_1}(a^*a+aa^*)$$ where $Q_1$ is the set of arrows of $Q$.
We denote by $e_i$ the idempotent of $\Lambda$ associated with the vertex $i$. We then have a natural functor 
$$\begin{array}{rcl}{\sf proj} \Lambda& \longrightarrow & \Ii ^\Pi (kQ) \\
e_i\Lambda&\mapsto& \prod_{p\geq 0}\tau_\Dd^pI_i\end{array}$$
where $\Ii^\Pi(kQ)$ is the closure of $\Ii(kQ)$ under countable products.
Composing this functor with the natural extension of $P$ to $\Ii^\Pi(kQ)$, we get a functor:

$$\begin{array}{rcl}{\sf proj} \Lambda&\longrightarrow & \Mm \\
e_i\Lambda &\mapsto& \bigoplus_{p\geq 0}\tau_B^p H_i.\end{array}$$

Therefore the restriction along this functor yields a functor $F:\mod \Mm\rightarrow
\mod\Lambda$. Moreover, since $\Mm$ has finitely many indecomposables, the functor $F$ takes its values in the full subcategory $\rm{f.l.}\Lambda$ formed by the $\Lambda$-modules of finite length.

This is an exact functor since it is a restriction. If $M$ is an $\Mm$-module, then the vector space $F(M)e_j$ is isomorphic to
$\bigoplus_{p\geq 0}M(\tau_B^pH_j)$. For $X$ in $\Mmm$, there exists $i\in Q_0$ and $q\geq 0$ such that $\tau^qH_i=X$. It is then easy to check that the image $F(S_X)$ of the simple associated to $X$ is the simple $\Lambda$-module $S_i$.

\subsubsection*{Fundamental propositions}

\begin{prop}\label{lemmefondamental}
For $X$ in $\Mmm$, there exists a functorial sequence in $\mod
\Lambda$ of the form
$$\xymatrix{0\ar[r] & F\circ i_*(X^\wedge)\ar[r] & F(H_0^\wedge)\ar[r] &
  F(H_1^\wedge)\ar[r] & F\circ i_*(X^\vee)\ar[r] & 0}$$
where $i_*:\mod \Mmm\rightarrow \mod \Mm$ is the right exact functor induced by $i:\Mmm\rightarrow \Mm$, and where $H_0$ and $H_1$ are in $add(H)$. 
\end{prop}

\begin{proof}
 Let $X$ be in $\Mmm$, and  $iX$ its image in $\Mm$. By $(1)$ of proposition \ref{propbase}, there exists a triangle functorial in $X$:
$$\xymatrix{iX\ar[r] & H_0\ar[r] & H_1\ar[r] & (iX)[1]}$$ with $H_0$ and $H_1$ in $add(H)$.
It yields a long exact sequence in $\mod \Mm$:
$$\xymatrix{0\ar[r] & (iX)^\wedge \ar[r] & H_0^\wedge\ar[r] & H_1^\wedge\ar[r] & \Ext^1_B(?,iX)_{|_\Mm}\ar[r] & \Ext^1_B(?,H_0)_{|_\Mm}\ar[r] & \cdots}.$$
By definition, the functor $\Ext^1_B(?,H_0)_{|_\Mm}$ is zero. The Auslander-Reiten formula gives us an isomorphism
$$\Ext^1_B(?,iX)_{|_\Mm}\simeq D\Hom_B(\tau_B^{-1}iX,?)_{|_\Mm}/\proj B.$$
Since $F$ is an exact functor, we get the following exact sequence in $f.l.\Lambda$:
$$\xymatrix{0\ar[r] & F((iX)^\wedge) \ar[r] & F(H_0^\wedge)\ar[r] & F(H_1^\wedge)\ar[r] & F((\tau_B^{-1}iX)^\vee/\proj B)\ar[r] & 0}$$
By definition, we have $F((iX)^\wedge)\simeq (F\circ i_*)(X^\wedge)$. For $j=1,\cdots, n$, we have an isomorphism:
$$
F((\tau_B^{-1}iX)^\vee/\proj B)e_j  \simeq \bigoplus_{p\geq 0} D\Hom_B(\tau_B^{-1}iX,\tau_B^p H_j)/\proj B.$$
For $p\geq 0$, we have $\tau^p_B(H_j)=\tau^{-1}_B(\tau_B^{p+1}H_j)$ if and only if $\tau_B^pH_j$ is not projective. Thus we have a vector space isomorphism
$$F((\tau_B^{-1}iX)^\vee/\proj B)e_j \simeq \bigoplus _{p\geq 0}D\Hom_B(\tau_B^{-1}iX,\tau_B^{-1}\tau_B^{p+1} H_j)/\proj B.$$
A morphism $f:\tau^{-1}X\rightarrow \tau^{-1}Y$ factorizes through a projective object if and only if $\tau(f):X\rightarrow Y$ is not zero. Thus we have:
\begin{eqnarray*}
F((\tau_B^{-1}iX)^\vee/\proj B)e_j & \simeq & \bigoplus_{p\geq 1} D\Hom_B(iX,\tau_B^p H_j)\\ 
&\simeq & \bigoplus_{p\geq 0} D\Hom_B(X,\tau_B^p H_j)/[add(H)]\\
&\simeq &(F\circ p^*)(X^\vee)e_j\simeq(F\circ i_*)(X^\vee)e_j.
\end{eqnarray*}

 Therefore we get this exact sequence in $f.l.\Lambda$, functorial in $X$:
 $$\xymatrix{0\ar[r] & (F\circ i_*)(X^\wedge) \ar[r] & F(H_0^\wedge)\ar[r] & F(H_1^\wedge)\ar[r] & (F\circ i_*)(X^\vee)\ar[r] & 0}$$
\end{proof}

\begin{prop}\label{propfond}
 Let $U$ and $V$ be indecomposables in $\Mmm$.  We have an isomorphism $$\Hom_{\Cc_A}(U^\wedge,V^\wedge)\simeq \bigoplus_{p\geq 0} \Mm(\tau_B^pU,V)/[add\tau_B^pH]$$ where $\Mm(\tau_B^pU,V)/[add\tau_B^pH]$ is the cokernel of the composition map $$\Mm(\tau_B^pU,\tau^p_BH)\otimes\Mm(\tau_B^pH,V)\longrightarrow \Mm(\tau_B^pU,V).$$
\end{prop}

We first show the following lemma:

\begin{lema}
 Let $e_U$ and $e_V$ be the idempotents of $A$ associated to the indecomposables $U$ and $V$. We have an isomorphism
$$e_U\Ext^2_A(DA,A)e_V\simeq \Mm(\tau_BU,V)/[add\tau_BH]$$ where $\Mm(\tau_BU,V)/[add\tau_BH]$ is the cokernel of the composition map $$\Mm(\tau_BU,\tau_BH)\otimes\Mm(\tau_BH,V)\longrightarrow \Mm(\tau_BU,V).$$
\end{lema}

\begin{proof}
We have the following isomorphisms:
\begin{eqnarray*} e_U\Ext^2_A(DA,A)e_V &=& \Ext^2_A(D(e_UA),Ae_V)\\
 & \simeq &\Hom_{\Dd(\Mmm)}(D\Mmm(U,?),\Mmm(?,V)[2]).
\end{eqnarray*}
Denote by $\underline{\Mm}$ the category $\Mm/\proj B$. The functor $\tau_B$ induces an equivalence of $k$-linear categories $\tau_B:\underline{\Mm}\rightarrow \Mmm.$ Thus we get the following isomorphisms
\begin{eqnarray*}
\Hom_{\Dd(\Mmm)}(D\Mmm(U,?),\Mmm(?,V)[2]) & \simeq & \Hom_{\Dd(\Mmm)}(D\underline{\Mm}(\tau_B^{-1}U,\tau_B^{-1}?),\underline{\Mm}(\tau_B^{-1}?,\tau_B^{-1}V)[2])\\&\simeq &\Hom_{\Dd(\underline{\Mm})}(D\underline{\Mm}(\tau_B^{-1}U,?),\underline{\Mm}(?,\tau_B^{-1}V)[2])\\
& \simeq & \Hom_{\Dd(\Mm)}(D\Mm(\tau_B^{-1}U,?)/\proj B,\Mm(?,\tau_B^{-1}V)/\proj B[2]) 
\end{eqnarray*}
But by the previous lemma, we know a projective resolution in $\mod \Mm$ of the module $D\Mm(\tau_B^{-1}U,?)/\proj B$. Namely, there exists an exact sequence in $\mod \Mm$ of the form:
$$\xymatrix{0\ar[r]& \Mm(?,U)\ar[r] & \Mm(?,H_0)\ar[r] &\Mm(?,H_1)\ar[r] & D\Mm(\tau_B^{-1}U,?)/\proj B\ar[r] & 0}$$ where $H_0$ and $H_1$ are in $add(H)$. Thus we get (using Yoneda's lemma)
\begin{eqnarray*}
\Hom_{\Dd(\Mmm)}(D\Mmm(U,?),\Mmm(?,V)[2]) & \simeq & \Hom_{\Dd(\Mm)}(\Mm(?,U),\Mm(?,\tau_B^{-1}V)/\proj B)/[add\Mm(?,H)]\\
&\simeq & \underline{\Mm}(U,\tau_B^{-1}V)/[add H]\\
&\simeq & \Mmm(\tau_B U,V)/[add \tau_B H].
\end{eqnarray*}
Since $V$ is in $\Mmm$, a non zero morphism of $\Mm(\tau_B U,V)$ cannot factorize through $add(H)$.
Thus we get $\Mmm(\tau_B U,V)/[add \tau_B H]\simeq \Mm(\tau_B U,V)/[add \tau_B H]$.

\end{proof}

\begin{proof}\textit{(of proposition \ref{propfond})}
In this proof, for simplicity we denote $\tau_B$ by $\tau$.
Let $\tilde{A}$ be the algebra $\End_{\Cc_A}(A)$. By proposition \ref{endomorphisme}, we have a vector space isomorphism
$$e_U\tilde{A}e_V\simeq e_UAe_V \oplus e_U\Ext^2_A(DA,A)e_V\oplus  e_U\Ext^2_A(DA,A)^{\otimes_A2}e_V\oplus \ldots$$
We prove by induction that $$e_U\Ext^2_A(DA,A)^{\otimes_Ap}e_V\simeq \Mm(\tau^pU,V)/[add\tau^pH].$$
For $p=0$, $e_UAe_V$ is isomorphic to $\Mmm(U,V)$ by Yoneda's lemma, and so to $\Mm(U,V)/[add(H)]$.
Suppose the proposition holds for an integer $p-1\geq 0$.
We then have $$e_u\Ext^2_A(DA,A)^{\otimes_Ap}e_V\simeq\sum_{W\in\ind(\Mmm)}e_u\Ext^2_A(DA,A)^{\otimes_Ap-1}e_W\ten e_W\Ext^2_A(DA,A)e_V.$$ The sum means here the direct sum modulo the mesh relations of the category $\Mmm$.
Thus this vector space is the sum over the indecomposables $W$ of $\Mmm$ of $$\Mm(\tau^{p-1}U,W)/[add(\tau^{p-1}H)]\otimes \Mm(\tau W,V)/[add(\tau H)]$$ modulo the mesh relations of $\Mmm$. This is isomorphic to the cokernel of the map $\varphi^{p-1}_{\tau^{p-1}U,W}\ten 1_{\tau W,V}+1_{\tau^{p-1}U,W}\ten \varphi^{1}_{\tau W,V}$ where 
$$\varphi^{j}_{X,Y}:\Mm(X,\tau^{j}H)\ten\Mm(\tau^{j}H,Y)\longrightarrow
 \Mm(X,Y)$$ is the composition map and where $$1_{X,Y}:\Mm(X,Y)\longrightarrow\Mm(X,Y)$$ is the identity.
The cokernel of this map is isomorphic to the cokernel of the map $\varphi^{p}_{\tau^p U,\tau W}\ten 1_{\tau W,V}+1_{U,\tau W}\ten \varphi^1_{\tau W,V}$. But we have an isomorphism 
$$\sum_{W\in\ind\Mmm}\Mm(\tau^pU,\tau W)\ten\Mm(\tau W,V)\simeq \Mm(\tau^pU,V).$$
Finally we get $$\Coker\left(\sum_{W\in\ind\Mmm}\varphi^{p}_{\tau^p U,\tau W}\ten 1_{\tau W,V}+1_{U,\tau W}\ten \varphi^{1}_{\tau W,V}\right)\simeq \Coker (\varphi_{\tau^pU,V}^{p}+\varphi^1_{\tau^pU,V}).$$
Furthermore, a morphism in $\Mm(\tau^pU,V)$ which factorizes through $\tau H$ factorizes through $\tau^p H$ since $H$ is a slice and $U$ is in $\Mmm$. Thus this cokernel is in fact isomorphic to the cokernel of $\varphi_{\tau^pU,V}^{p}$ that is to say to the space 
$$\Mm(\tau^pU,V)/[add\tau^pH].$$
\end{proof}

\subsection{Case where $B$ is hereditary}

\subsubsection*{Results of Geiss, Leclerc and Schr{\"o}er } 

Let $Q$ be a finite connected quiver without oriented cycles with
 $n$ vertices. Denote by $\Pp$ the postprojective component of the
 Auslander-Reiten quiver of $\mod kQ$, and by $P_1, \ldots , P_n$ the indecomposable
 projectives. 

\begin{dfa}[Geiss-Leclerc-Schr{\"o}er, \cite{Gei4}]
A $kQ$-module $M=M_1\oplus\cdots\oplus M_r$, where the $M_i$ are pairwise non isomorphic indecomposables, is called \emph{initial}
if the following conditions hold:
\begin{itemize}
\item for all $i=1,\ldots,r$, $M_i$ is postprojective;
\item if $X$ is an indecomposable $kQ$-module with $\Hom_{kQ}(X,M)\neq
  0$, then $X$ is in $add(M)$;
\item and $P_i\in add(M)$ for each indecomposable projective $kQ$-module
  $P_i$.
\end{itemize}
We define the integers $t_i$ as $$t_i=max \{ j\geq 0 | \tau^{-j}(P_i)\in add(M)-\{0\}\}.$$
\end{dfa}
   
Denote by $\Lambda$ the preprojective algebra associated to
$Q$. There is a canonical embedding of algebras 
$\xymatrix{kQ\ar@{^(->}[r] & \Lambda}.$ Denote by $\pi_Q:\mod
\Lambda\rightarrow \mod kQ$ the corresponding restriction functor.

\begin{thma}[Geiss-Leclerc-Schr{\"o}er, \cite{Gei4}]
Let $M$ be an initial $kQ$-module, and let $\Cc_M=\pi_Q^{-1}(add(M))$
be the subcategory of all $\Lambda$-modules $X$ with $\pi_Q(X)\in
add(M)$.
 The following holds:

\hspace{.2in}(i) the category $\Cc_M$ is a Frobenius category with $n$
projective-injectives;

\hspace{.2in}(ii) the stable category $\underline{\Cc}_M$ is a
$2$-CY triangulated category.
\end{thma}

Recall from Ringel \cite{Rin3} that the category $\mod \Lambda$ can be seen as
$\mod kQ (\tau^{-1}, 1)$. The objects are pairs $(X,f)$ where $X$ is
in $\mod kQ$ and  $f:\tau^{-1}X\rightarrow
X$ is a morphism in $\mod kQ$. The morphisms $\varphi$ between $(X,f)$
and $(Y,g)$ are commutative squares:
$$\xymatrix{\tau^{-1} X\ar[r]^f \ar[d]_{\tau^{-1}\varphi}&
  X\ar[d]^\varphi \\ \tau^{-1} Y \ar[r]^g & Y}$$
The image of an object $(X,f)$ under $\pi_Q:\mod \Lambda\rightarrow \mod kQ$
is then the module $X$.

Let $X=\tau ^{-l}P_i$ be an indecomposable summand of an initial
module $M$.  Let $R_X=(Y,f)$ be the following object in $\mod kQ
(\tau^{-1}, 1)\simeq \mod \Lambda$:

$$Y=\bigoplus_{j=0}^l \tau^{-j}P_i \quad {\rm and} \quad f:\bigoplus_{j=1}^{l+1} \tau^{-j}P_i \longrightarrow 
  \bigoplus_{j=0}^l \tau^{-j}P_i$$ is given by the matrix 
$$f=\left( \bsm 0 &  & & \\ 1& \ddots & &  \\ &\ddots
    & \ddots &  \\ & &1 & 0\esm\right).$$

\begin{prop}[Geiss-Leclerc-Schr{\"o}er,\cite{Gei4}]
The category $\Cc_M$ has a canonical maximal rigid object
$R=\bigoplus_{X\in \ind add(M)} R_X$. The projective-injectives of
$\Cc_M$ are the $R_{\tau^{-t_i}P_i}$, $i=1,\ldots ,n$. Therefore, $R$ is a cluster-tilting object in $\underline{\Cc}_M$.
\end{prop}

\subsubsection*{Endomorphism algebra of the cluster-tilting object}
Let $Q$ be a connected quiver without oriented cycles and denote by
$B$ the path algebra $kQ$. Let $M$ be
an initial $B$-module. Let $H$ be the following postprojective slice 
$H=\bigoplus_{i=1}^n \tau^{-t_i}P_i$ of $\mod B$.
Let $Q'$ be the quiver such that $\End_B(H)$ is isomorphic to $kQ'$.

Let us define, as in the previous section, the subcategory $\Mm$ of
$\Dd^b(\mod kQ)$ as $$\Mm
=\{X\in \mod kQ \ /  \Ext^1_{B}(X,H)=0\}.$$
It is then obvious that $\Mm=add(M)$.
As previously, we denote by $\Lambda$ the preprojective algebra associated with $Q'$. It is isomorphic to the one associated with $Q$ because $Q$ and
$Q'$ have the same underlying graph. Recall that we have
$\Mmm=\Mm/add(H)$, and that $A=\End_B(\Mmm)$ is an algebra of global
dimension $2$. Note that in this case $\tau_B$ and $\tau_\Dd$
coincide on the objects of $\mod B$ which have no projective direct summands since $B$ is hereditary. We will denote it by $\tau$ in this section.

\begin{lema}
 Let $U$ and $V$ be indecomposables in $\Mmm$. We have $$\Hom_\Lambda(R_U,R_V)\simeq \bigoplus_{j\geq 0} \Mm(\tau^jU,V).$$
\end{lema}

\begin{proof}
 Let $P$ and $Q$ be projective indecomposables such that $U=\tau^{-q}Q$ and $V=\tau^{-p}P$. 

\textit{Case 1: $p\leq q$}\\
An easy computation gives the following equalities 
\begin{eqnarray*}\Hom_\Lambda(R_U,R_V) & \simeq & \bigoplus_{j=0}^p\Mm(Q,\tau^{-j}P) 
  \simeq  \bigoplus_{j=0}^p\Mm(\tau^{-p+j}Q,\tau^{-p}P)\\
  & \simeq & \bigoplus_{j=0}^p\Mm(\tau^{-p+j+q}(\tau^{-q}Q),\tau^{-p}P)  \simeq  \bigoplus_{j=q-p}^q\Mm(\tau^{j}U,V).
\end{eqnarray*}
Since $\Mm(\tau^kU,V)$ vanishes for $k\leq q-p+1$ and since $\tau^kU$ vanishes for $k\geq q+1$  we get an isomorphism
$$ \Hom_\Lambda(R_U,R_V)\simeq\bigoplus_{j\geq 0}\Mm(\tau^{j}U,V).$$

\textit{Case 2: $p>q$}\\
In this case, a morphism from $R_U$ to $R_V$ is given by morphisms $a_j\in \Mm(Q,\tau^{-j}P)$, with $j=0,\ldots, p$ such that $\tau^{-q+1}a_j=0$ for $j=0,\ldots,p-q-1$. But since $\tau^{-q+1-j}P$ is not zero for $j=0,\ldots,p-q-1$, the morphism $\tau^{-q+1}a_j:\tau^{-q+1}Q\rightarrow \tau^{-q+1-j}P$ vanishes if and only if $a_j$ vanishes. Thus we get \begin{eqnarray*}\Hom_\Lambda(R_U,R_V) & \simeq  & \bigoplus_{j=p-q}^p\Mm(Q,\tau^{-j}P) \simeq  \bigoplus_{j=p-q}^p\Mm(\tau^{-p+j}Q,\tau^{-p}P)\\
  & \simeq & \bigoplus_{j=p-q}^p\Mm(\tau^{-p+j+q}(\tau^{-q}Q),\tau^{-p}P)   \simeq  \bigoplus_{j=0}^q\Mm(\tau^{j}U,V).\end{eqnarray*}
Since $\tau^jU$ vanishes for $j\geq q+1$ we get $$ \Hom_\Lambda(R_U,R_V)\simeq\bigoplus_{j\geq 0}\Mm(\tau^{j}U,V).$$

\end{proof}

\begin{cora}\label{addTaddA}
 Let $U$ and $V$ be indecomposable objects in $\Mmm$. We have $$\Hom_{\underline{\Cc}_M}(R_U,R_V)\simeq e_U\tilde{A}e_V$$ and therefore the algebras $\tilde{A}$ and $\End_{\underline{\Cc}_M}(R)$ are isomorphic.
\end{cora}

\begin{proof}
 The projective-injectives in the category $\Cc_M$ are the $R_{H_i}$ with $i=1,\ldots, n$. Denote by $R_H$ the sum $\bigoplus_{i=1}^nR_{H_i}$. Thus $\Hom_{\underline{\Cc}_M}(R_U,R_V)$ is the cokernel of the composition map $$\Hom_{\Cc_M}(R_U,R_H)\otimes\Hom_{\Cc_M}(R_H,R_V)\longrightarrow\Hom_{\Cc_M}(R_U,R_V).$$
By the previous lemma  this map is isomorphic to the following 
$$\xymatrix{\bigoplus_{i,j\geq 0}\Mm(\tau^iU,H)\otimes\Mm(\tau^jH,V)\ar[rr]^-\Phi &&\bigoplus_{p\geq 0}\Mm(\tau^pU,V)}$$
Given two morphisms $f\in\Mm(\tau^iU,H)$ and $\Mm(\tau^jH,V)$, $\Phi(f\otimes g)$ is the composition $\tau^jf\circ g\in\Mm(\tau^{i+j}U,V)$. Thus the cokernel of this map is the cokernel of the map $$\xymatrix{\bigoplus_{p\geq 0}\bigoplus_{i=0}^p\Mm(\tau^pU,\tau^iH)\otimes\Mm(\tau^iH,V)\ar[rr]^-\Phi &&\bigoplus_{p\geq 0}\Mm(\tau^pU,V)}. $$
Since $H$ is a slice and since $U$ is in $\Mmm$, a morphism in $\Mm(\tau^pU,V)$ which factorizes through $\tau^iH$ with $i\leq p$ factorizes through $\tau^pH$. Finally we get $$\Hom_{\underline{\Cc}_M}(R_U,R_V)\simeq \bigoplus_{p\geq 0}\Mm(\tau^pU,V)/[add\tau^pH],$$ and we conclude using proposition \ref{propfond}.
\end{proof}

\subsubsection*{Triangle equivalence}

\begin{thma}\label{casGLS}
The functor $F\circ i_*:\mod \Mmm\rightarrow f.l.\Lambda$ yields a
triangle equivalence between $\Cc_{\Mmm}$ and $\underline{\Cc}_M$.
\end{thma}

\begin{proof}
Let $X=\tau_B^{-l}P_i$ be an indecomposable of $\Mm$. Let $X^\wedge$ be the projective $\Mm$-module
$\Hom_B(?,X)_{|_\Mm}$. The underlying vector space of $F(X^\wedge)$ is 
\begin{eqnarray*}F(X^\wedge)
 &\simeq& \bigoplus_{q\geq 0}\Hom_B(\tau_B^q H,\tau_B^{-l}P_i)\simeq \bigoplus_{q\geq
   0}\Hom_B(\tau_B^{-q}B,\tau_B^{-l}P_i)\\
 &\simeq&\bigoplus_{q\geq 0}\Hom_B(B,\tau_B^{q-l}P_i)\simeq\bigoplus_{q=0}^l \tau_B^{-q}P_i.
\end{eqnarray*} 

It is then not hard to see that $F(X^\wedge)$ is equal to $R_X$.  
Thus each projective $X^\wedge$ is sent onto an object of $\Cc_M$. Therefore
$F$ induces a functor $F:\Dd^b(\Mm)\rightarrow \Dd^b(\Cc_M)$. Moreover for $i=1,\ldots,n$, $F(H_i^\wedge)$ is equal to $R_{\tau^{-t_i}P_i}$,\emph{i.e.} a projective-injective of $\Cc_M$.  
We have the following composition:
$$\xymatrix{\Dd^b(\Mmm)\simeq\Dd^b( A)\ar@(dl,dr)_{?\lten_A
    DA[-2]}\ar[r]^-{i_*} & \Dd^b(\Mm)\ar[r]^-F & \Dd^b(\Cc_M)\ar[r]^-{\pi} &
  \Dd^b(\Cc_M)/\per\Cc_M\simeq \underline{\Cc}_M}$$
The functor $F\circ i_*$ is clearly isomorphic to the
left derived tensor product with the $A$-$\Lambda$-bimodule
$R=F\circ i_*(A)$. By proposition \ref{lemmefondamental}, for $X$ in $\overline{\mathcal{M}}$,
we have the following exact sequence, functorial in $X$:
$$\xymatrix{0\ar[r] & F\circ i_*(X^\wedge)\ar[r] & F(H_0^\wedge)\ar[r] &
  F(H_1^\wedge)\ar[r] & F\circ i_*(X^\vee)\ar[r] & 0}$$ with $H_0$ and
$H_1$ in $add(H)$. It yields a morphism
$$F\circ i_* (DA)\rightarrow  F\circ i_* (A)[2]$$
in the derived category of $A$-$\Lambda$-bimodules. Since the
objects $F(H_0^\wedge)$ and $F(H_1^\wedge)$ vanish in the
stable category $\underline{\mathcal{C}}_M$, the image $$F\circ i_* (DA)\rightarrow  F\circ i_*(A)[2]$$ of this morphism in the category of $A$-$\mathcal{B}$-bimodules
is invertible, where $\mathcal{B}$ is a dg category whose perfect
derived category is algebraically equivalent to the stable
category $\underline{\mathcal{C}}_M$. In other words, in
the derived category $\mathcal{D}(A^{op}\otimes \mathcal{B})$,
we have an isomorphism
$$DA\lten_A \pi F i_*(A) \simeq \pi F i_*(A)[-2].$$
By the universal property of the orbit category, we have the factorization

$$\xymatrix@R=-.3pc{\Dd^b( \Mmm)\ar[rr]^{?\lten_A \underline{R}}\ar[dr]
  && \underline{\Cc}_M. \\ & \Cc_{\Mmm} \ar@{..>}[ur] &}$$

This factorization is an algebraic functor between $2$-CY
categories which sends the cluster-tilting object $A$ onto the cluster-tilting object $\underline{R}$. Moreover by corollary \ref{addTaddA}, it yields an equivalence between the categories $add(A)$ and $add(\underline{R})$. Thus it is an algebraic triangle
equivalence.

 \end{proof}

Note that if $M$ is the initial module $kQ\oplus \tau^{-1} kQ$, Geiss, Leclerc and Schr{\"o}er proved, using a result of Keller and Reiten \cite{Kel4}, that the $2$-CY category $\underline{\Cc}_M$ is triangle equivalent to the cluster category $\Cc_Q$. Here, $H$ is $\tau^{-1}kQ$ and then $\Mmm$ is $kQ$, so we get another proof of this fact. 

\subsection{Relation with categories $\Sub\Lambda/\Ii_w$}

\subsubsection*{Results of Buan, Iyama, Reiten and Scott} 

Let $Q$ be a finite connected quiver without oriented cycles and $\Lambda$ the associated preprojective algebra. We denote by $\{1,\ldots, n\}$ the set of vertices of $Q$. For a vertex $i$ of $Q$, we denote by $\Ii_i$ the ideal $\Lambda(1-e_i)\Lambda$ of $\Lambda$. We denote by $W$ the \emph{Coxeter group} associated to the quiver $Q$. The group $W$ is defined by the generators ${1,\ldots,n}$ and the relations:
\begin{itemize}
\item $i^2=1$ for all $i$ in $\{1,\ldots, n\}$;
 \item $ij=ji$ if there are no arrows between the vertices $i$ and $j$;
\item $iji=jij$ if there is exactly one arrow between $i$ and $j$. 
\end{itemize}

Let $w=i_1i_2\ldots i_r$ be a $W$-reduced word. For $m\leq r$, let $\Ii_{w_m}$ be the following ideal:
$$\Ii_{w_m}=\Ii_{i_m}\ldots \Ii_{i_2}\Ii_{i_1}.$$
For simplicity we will denote $\Ii_{w_r}$ by $\Ii_w$. The category $\Sub \Lambda/\Ii_w$ is the subcategory of $\rm{f.l.}\Lambda$ generated by the sub-$\Lambda$-modules of $\Lambda/\Ii_w$.

\begin{thma}[Buan-Iyama-Reiten-Scott \cite{Bua2}]
The category $\Sub \Lambda/\Ii_w$ is a Frobenius category and its stable category $\underline{\Sub} \Lambda/\Ii_w$ is $2$-CY. The object $T_w=\bigoplus_{m=1}^r e_{i_m}\Lambda/\Ii_{w_m}$ is a cluster-tilting object.
 \end{thma}
Note that this theorem is written only for non Dynkin quivers in \cite{Bua2}, but the Dynkin case is an easy consequence of theorem II.2.8 and corollary II.3.5 of \cite{Bua2}. 

\subsubsection*{Construction of a reduced word}

Let $B$ be a concealed algebra, and $H$ a postprojective slice in $\mod B$. Let $Q$ the quiver of $\End_B(H)$. It is a finite quiver without oriented cycles. We denote by $\{1,\ldots ,n\}$ its set of vertices and by  
 $\Lambda$ its preprojective algebra. We define as previously  
$\Mm
=\{X\in \mod B \ /  \Ext^1_{B}(X,H)=0\}$.

Let us order the indecomposables $X_1,\ldots, X_N$ of $\Mm$ in such a way: if the morphism space $\Hom_B(X_i,X_j)$ does not vanish, $i$ is smaller than $j$. This is possible since $Q$ has no oriented cycles.

By proposition \ref{propbase}, for $X_i\in \Mm$ there exists a unique $q\geq 0$ such that $\tau_B^{-q}X_i\simeq H_{\varphi(i)}$ for a certain integer $\varphi(i)$. So we get a function $\varphi:\{1,\ldots,N\}\rightarrow \{1,\ldots, n\}$. Let $w$ be the word $\varphi(1)\varphi(2)\ldots\varphi(N)$.

\begin{prop}\label{Wreduit}
The word $w$ is $W$-reduced.
\end{prop}

\begin{proof}
The proof is in several steps:\medskip

\textit{Step 1: For two integers $i<j$ in $\{1,\ldots,N\}$, we have $\varphi(i)=\varphi(j)$ if and only if there exists a positive integer $p$ such that $X_i=\tau_B^pX_j$.}\medskip

\textit{Step 2: The element w of the Coxeter group does not depend on the order on the indecomposables of $\Mm$.}\smallskip\\
Let $i$ be in $\{1,\ldots,N-1\}$. Assume there is an arrow $\varphi(i)\rightarrow\varphi(i+1)$ in $Q$. We show that there is an arrow $X_i\rightarrow X_{i+1}$  in the Auslander-Reiten quiver of $\Mm$.
By proposition \ref{propbase}, there exist positive integers $p$ and $q$ such that $X_i=\tau_B^qH_{\varphi(i)}$ and $X_{i+1}=\tau_B^pH_{\varphi(i+1)}$. By hypothesis there is an arrow between $H_{\varphi(i)}$ and $H_{\varphi(i+1)}$. Thus we want to show that $p$ is equal to $q$.

 Suppose that $p\geq q+1$, then since $\Mm$ is closed under $\tau_B$, the objects $\tau^q_BH_{\varphi(i+1)}$ and $\tau^{q+1}_BH_{\varphi(i+1)}$ are non zero and are in $\Mm$. Let $l$ be the integer in $\{1,\ldots,N\}$ such that $X_l=\tau^{q+1}_BH_{\varphi(i+1)}$. We have an arrow $$X_i=\tau_B^qH_{\varphi(i)}\rightarrow \tau_B^qH_{\varphi(i+1)}=\tau_B^{-1}X_l.$$ Thus, by the property of the AR-translation, there is an arrow $X_l\rightarrow X_i.$ Thus $i$ should be strictly greater than $l$. But by step 1, and the hypothesis $p\geq q+1$, we have $i+1\leq l$. This is a contradiction.

The cases $q\geq p+1$, and  $\varphi(i+1)\rightarrow\varphi(i)$ in $Q$ can be solved in the same way.\medskip

\textit{Step 3: It is not possible to have $\varphi(i)=\varphi(i+1)$.}\smallskip\\
Suppose we have $\varphi(i)=\varphi(i+1)$. By step 1 there exists a positive integer $p$ such that $X_i=\tau_B^pX_{i+1}$. Suppose that $p$ is  $\geq 2$, then $\tau_BX_{i+1}=\tau_B^{-p+1}X_i$ is in $\Mm$, it is isomorphic to an $X_k$ for an integer $k$ with $\varphi(k)=\varphi(i)$. But $k$ must be strictly greater than $i$ and strictly smaller than $i+1$ which is clearly impossible. 
Thus $p$ is equal to 1. There should exist an $X_l$ in $\Mm$ such that $\Hom(X_i,X_l)\neq 0$ and $\Hom(X_l,X_{i+1})\neq 0$. Thus $l$ must be strictly between $i$ and $i+1$ which is impossible. \medskip

\textit{Step 4: It is not possible to have $\varphi(i)=\varphi(i+2)$ and $\varphi(i+1)=\varphi(i+3)$ with exactly one arrow in $Q$ between $\varphi(i)$ and $\varphi(i+1)$.}\smallskip\\
In this case we have, by step 1, $X_i=\tau_B^pX_{i+2}$ and $X_{i+1}=\tau_B^qX_{i+3}$. By the same argument as in step 3, $p$ and $q$ have to be equal to $1$. Thus the AR quiver of $\Mm$ has locally the following form:
$$
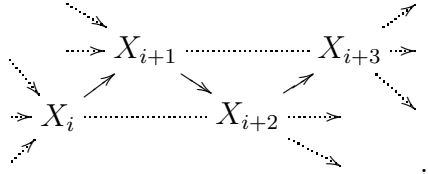
\xymatrix@-1.4pc{&\ar@{.>}[dr]&&&&\\ \ar@{.>}[dr]&\ar@{.>}[r]&X_{i+1}\ar@{.}[rr]\ar[dr]&&X_{i+3}\ar@{.>}[ur]\ar@{.>}[r]\ar@{.>}[dr]&\\
\ar@{.>}[r]& X_i\ar@{.}[rr]\ar[ur]&&X_{i+2}\ar[ur]\ar@{.>}[r]\ar@{.>}[dr]&&\\
\ar@{.>}[ur]&&&&&.}$$
The module $X_{i+1}$ is the unique direct predecessor of $X_{i+2}$. Indeed, suppose there is an $X_k$ with an arrow $X_k\rightarrow X_{i+2}$. Thus there is an arrow $\tau_BX_{i+2}=X_i\rightarrow X_k$ and $k$ must be strictly between $i$ and $i+2$. By the same argument, there is only one arrow with tail $X_{i+3}$, one arrow with source $X_i$ and one arrow with source $X_{i+1}$.
 Thus we have the following AR sequences in $\mod B$:
$$\xymatrix@-1.1pc{0\ar[r] & X_i\ar[r] & X_{i+1}\ar[r] & X_{i+2}\ar[r] & 0}\quad \textrm{and} \quad \xymatrix@-1.1pc{0\ar[r] & X_{i+1}\ar[r] & X_{i+2}\ar[r] & X_{i+3}\ar[r] & 0}$$
which is clearly impossible.\medskip

\textit{Step 5: There is no subsequence of type jkjlkl in w with an arrow between j and k and an arrow between k and l}\smallskip\\
Suppose we have $\varphi(i)=\varphi(i+2)=j$, $\varphi(i+1)=\varphi(i+4)=k$ and $\varphi(i+3)=\varphi(i+5)=l$. As previously, we have $X_i=\tau_BX_{i+2}$, $X_{i+1}=\tau_BX_{i+4}$ and $X_{i+3}=\tau_BX_{i+5}$. There is an arrow $X_{i+1}\rightarrow X_{i+2}$ so there is an arrow $X_{i+2}\rightarrow X_{i+4}$. There is an arrow $X_{i+3}\rightarrow X_{i+4}$ thus there is an arrow $X_{i+1}\rightarrow X_{i+3}$. As in step 4 it is easy to see that the AR quiver of $\Mm$ locally looks like:
$$
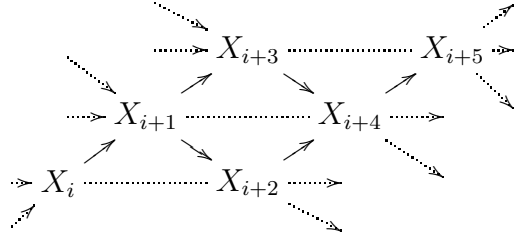
\xymatrix@-1.4pc{&&\ar@{.>}[dr]&&&&\\ &\ar@{.>}[dr]&\ar@{.>}[r]&X_{i+3}\ar@{.}[rr]\ar[dr]&&X_{i+5}\ar@{.>}[ur]\ar@{.>}[r]\ar@{.>}[dr]&\\
&\ar@{.>}[r]& X_{i+1}\ar@{.}[rr]\ar[ur]\ar[dr]&&X_{i+4}\ar[ur]\ar@{.>}[r]\ar@{.>}[dr]&&\\
\ar@{.>}[r]& X_{i}\ar@{.}[rr]\ar[ur]&&X_{i+2}\ar[ur]\ar@{.>}[r]\ar@{.>}[dr]&&&\\
\ar@{.>}[ur]&&&&&&\\
}$$
Thus we have the 3 following AR sequences in $\mod B$:
$$\xymatrix@-1.1pc{0\ar[r] & X_i\ar[r] & X_{i+1}\ar[r] & X_{i+2}\ar[r] & 0} \quad \xymatrix@-1.1pc{0\ar[r] & X_{i+3}\ar[r] & X_{i+4}\ar[r] & X_{i+5}\ar[r] & 0}$$
$$\textrm{and}\quad \xymatrix{0\ar[r] & X_{i+1}\ar[r] & X_{i+3}\oplus X_{i+2}\ar[r] & X_{i+4}\ar[r] & 0}$$
A simple argument of dimension permits us to conclude that $X_i$ and $X_{i+5}$ must be zero, that is a contradiction.

By the second step, we know that using the relation of commutativity is the same as changing the order on the indecomposables of $\Mm$. Moreover we just saw that locally we can not reduce the word $w$. Thus it is reduced.

\end{proof}

\subsubsection*{Image of the cluster-tilting object}

Let $F:\mod\Mm\rightarrow \textrm{f.l.}\Lambda$ be the functor constructed in section \ref{construction}.

\begin{prop}\label{imagedeF}
For $i=1,\ldots, N$, we have an isomorphism in $\textrm{f.l.}\Lambda$:
$$F(X_i^\wedge)\simeq e_{\varphi (i)} \Lambda/\Ii_{w_{i}}$$ where $w_i$ is the word $\varphi(1)\cdots \varphi(i)$. 
\end{prop}

\begin{proof}
The functor $F$ is right exact and sends the simple functor $S_{X_i}$ onto the simple $S_{\varphi(i)}$. Since $F(X_i^\wedge)$ surjects onto  $F(S_{X_i})$, there is a morphism $e_{\varphi(i)}\Lambda\rightarrow F(X_i^\wedge)$. Explicitly, we will take the morphism  given in this way:

The object $X_i$ is of the form $\tau_B^qH_{\varphi(i)}$ for a $q\geq 0$. 
If $j$ is in $\{1,\ldots,n\}$, the vector space $e_{\varphi(i)}\Lambda e_j$ is isomorphic to $\prod_{p\geq 0} \Hom_{kQ}(\tau^p_\Dd I_j, I_{\varphi(i)})$ where $I_j$ is the injective indecomposable module of $\mod kQ$ corresponding to the vertex $j$.
Let $f$ be a morphism in $\Hom_{kQ}(\tau^p_\Dd I_j, I_{\varphi(i)})$, then $\tau_\Dd^q(f)$ is a morphism in 
$\Hom_{kQ}(\tau^{p+q}_\Dd I_j, \tau_\Dd^q I_{\varphi(i)})$, and then $P(\tau_\Dd^q f)=\tau_B^q P(f)$ is a morphism in $\Mm$ from $\tau_B^{p+q} H_j$ to $\tau_B^q H_{\varphi(i)}=X_i$, thus is in $F(X_i^\wedge)e_j$. \medskip\\

\textit{Step 1: The morphism $e_{\varphi (i)}\Lambda \rightarrow F(X_i^\wedge)$ vanishes on the ideal $\Ii_{w_i}$.}\smallskip

A word $j_1j_2\cdots j_r$ will be called a \emph{subword} of $w_i$ if there exist integers $1\leq l_1<l_2<\cdots <l_r\leq i$ such that $j_1j_2\cdots j_r=\varphi(l_1)\varphi(l_2)\cdots\varphi(l_r)$.
It is easy to check that the vector space $e_{\varphi(i)}\Ii_{w_i}e_j$ is generated by the paths from $j$ to $\varphi(i)$ such that there exists a factorization $$\xymatrix{j\ar@{~>}[r] &j_1\ar@{~>}[r] & j_2\ar@{~>}[r] & \cdots \ar@{~>}[r] & j_r\ar@{~>}[r] & \varphi(i)}$$ with $jj_1j_2\cdots j_r\varphi(i)$ not a subword of $w_i$. 

Let $f$ be a morphism $\tau_\Dd^pI_j\rightarrow I_{\varphi(i)}$ in $\Ii(kQ)$ given by such a path. Assume that the image  $P(\tau_\Dd^q f)$ of $f$ in $F(X_i^\wedge)$ is non zero. Let $$\xymatrix{\tau^p_\Dd I_j\ar[r]^{f_0} & \tau_\Dd^{p_1}I_{j_1}\ar[r]^{f_1} &\tau_\Dd^{p_2}I_{j_2}\ar[r]^{f_2}&\cdots \ar[r] &\tau_\Dd^{p_r}I_{j_r}\ar[r]^{f_r} & I_{\varphi(i)}}$$ be the factorization of $f$ given by the above  factorization of the path. Therefore $P(\tau_\Dd^q f)$ is equal to the composition
$$\xymatrix{\tau^{p+q}_B H_j\ar[r] & \tau_B^{p_1+q}H_{j_1}\ar[r] &\tau_B^{p_2+q}H_{j_2}\ar[r]&\cdots \ar[r] &\tau_B^{p_r+q}H_{j_r}\ar[r] & \tau_B^q H_{\varphi(i)}=X_i}.$$
Since $P(\tau_\Dd^q f)$ is not zero, all morphisms $P(\tau_\Dd^q f_l)$ are not zero, and all objects $\tau_B^{p_l+q}H_{j_l}$ are non zero. Thus the objects $\tau_B^{p_l+q}H_{j_l}$ are of the form $X_{h_l}$ with $h_0<h_1<\cdots <h_r<i$. Furthermore, we have $\varphi(h_l)=j_l$. Thus $j j_1\cdots j_r \varphi(i)=\varphi(h_0)\varphi(h_1)\cdots\varphi(h_r)\varphi(i)$ is a subword of $w_i$. This contradiction shows that the image of $f$ in $F(X_i^\wedge)$ must be zero.\medskip\\

\textit{Step 2: The morphism $e_{\varphi(i)}\Lambda \rightarrow F(X_i^\wedge)$ is surjective.}\smallskip

Let $f$ be a morphism $\tau_B^{p+q}H_j\rightarrow \tau_B^q H_{\varphi(i)}=X_i$ in $\Mm$. Hence $\tau_B^{-q} f$ is a morphism $\tau_B^{p}H_j\rightarrow H_{\varphi(i)}$ in $\Mm$. Since $P$ is full (cf. prop.~\ref{constructionP}), there exists a morphism $g:\tau_\Dd^pI_i\rightarrow I_{\varphi(i)}$ such that $P(g)=\tau_B^{-q}f$.  Thus we have $P(\tau_\Dd^qg)=\tau_B^qP(g)=f$.
\medskip\\

\textit{Step 3: The morphism $e_{\varphi(i)}\Lambda/\Ii_{w_i} \rightarrow F(X_i^\wedge)$ is injective.}\smallskip

  Let $f$ be a non zero morphism $\tau_\Dd^pI_j\rightarrow I_{\varphi (i)}$ in $\Ii(kQ)$ such that $P(\tau_\Dd^qf)$ is zero.  By lemma \ref{kernel}, we can assume that there exists a factorization of $\tau_\Dd^qf$ of the form $$\xymatrix{\tau_\Dd^{q+p}I_j\ar[r]^h & Y\ar[r]^-g & \tau^q_\Dd I_{\varphi (i)}}$$ with $Y$ indecomposable and $P(Y)=0$. The object $Y$ is of the form $\tau_\Dd^{h}I_l$ with $h\geq q$ and we have $\tau_B^{h}H_l=0$.

The morphism $g$ is a sum of compositions of irreducible morphisms between indecomposables. Let  
$$\xymatrix{\tau_\Dd^{h}I_l\ar[r]^{g_0} & Y_1\ar[r]^{g_1} & Y_2\ar[r]^{g_2} &\cdots \ar[r]& Y_s\ar[r]^{g_s} & \tau^q_\Dd I_{\varphi(i)}}$$ be such a summand of $g$. The objects $Y_k$, $1\leq k\leq s$ are indecomposable and  so are of the form $\tau_\Dd^{r_k}I_{j_k}$, and the morphisms $g_k$, $0\leq k\leq s$ are irreducible. We will show that the word $lj_1j_2\ldots j_s\varphi(i)$ is not a subword of $w_i$. Without loss of generality, we may assume that for $1\leq k\leq s$, $P(Y_k)$ is not zero, so there exist integers $l_k$ such that $P(Y_k)=X_{l_k}$. Since the morphisms $g_k$ are irreducible, $P(g_k)$ does not vanish, and we have $1\leq l_1<l_2<\cdots<l_s<i$. The word $j_1j_2\ldots j_s\varphi(i)$ is equal to the word $\varphi(l_1)\varphi(l_2)\cdots\varphi(l_s)\varphi(i)$, so $j_1j_2\ldots j_s\varphi(i)$ is a subword of $w_i$.\medskip\\

\textit{Substep 1: The sequence $1\leq l_1<l_2<\cdots<l_s<i$ is the maximal element of the set $\{1\leq i_1<i_2<\cdots<i_s<i_{s+1}\leq i \quad | \quad \varphi(i_1)=j_1, \ldots, \varphi(i_s)=j_s, \varphi(i_{s+1})=\varphi(i)\}$ for the lexicographic order.}\smallskip\\
We prove by decreasing induction that $l_k$ is the maximal integer with $l_k<l_{k+1}$ and $\varphi(l_k)=j_k$. For $k=s+1$ it is obvious.
Now suppose there exists an integer $i_k$ such that $\varphi(l_k)=\varphi(i_k)=j_k$ and $l_k<i_k<l_{k+1}$.
Thus by step 1 of proposition \ref{Wreduit}, there exists an integer $r\geq 1$ such that $X_{l_k}=\tau_B^rX_{i_k}$. The morphism $P(g_k):X_{l_k}\rightarrow X_{l_{k+1}}$ is irreducible, so there exists a non zero irreducible morphism $X_{l_{k+1}}\rightarrow \tau_B^{-1}X_{l_k}$. The object $\tau_B^{-1}X_{l_k}$ is in $\Mm$ since $X_{l_{k}}$ and $\tau_B^{-r}X_{l_k}=X_{i_k}$ are in $\Mm$. It is of the form $X_t$, and we have $l_{k+1}<t$. Since $r$ is $\geq 1$, $t$ is $\leq i_k$ by step 1 of proposition ~\ref{Wreduit}. This implies $l_{k+1}<i_k$ which is a contradiction.\medskip\\

\textit{Substep 2: $l$ does not belong to the set $\{ \varphi(1),\varphi(2), \ldots, \varphi(l_1-1)\}$.}\smallskip\\
Suppose that there exists an integer $1\leq k\leq N$ such that $\varphi(k)$ is equal to $l$. Thus there exists an integer $r\geq 0$ such that $X_k$ is equal to $\tau_B^rH_l$. Since $\tau_B^{h}H_l=P(\tau_\Dd^hI_l)$ is zero, $r$ is $\leq h-1$. 

Since the morphism $g_0:\tau_\Dd^hI_l\rightarrow Y_1$ is an irreducible morphism of $\Ii(kQ)$, there exists an irreducible morphism $Y_1\rightarrow \tau_\Dd^{h-1} I_l$ in $\Ii(kQ)$. Thus there exists an irreducible morphism $\tau_\Dd^{r-h+1}Y_1\rightarrow \tau_\Dd^r I_l$ in $\Ii(kQ)$. The object $P(\tau_\Dd^rI_l)=\tau_B^rH_l=X_k$ is not zero and lies in $\Mm$, so the object $P(\tau_\Dd^{r-h+1}Y_1)=\tau_B^{r-h+1}X_{l_1}$ is not zero and lies in $\Mm$ since $\Mm$ is stable by kernel. Thus there is an irreducible morphism $\tau_B^{r-h+1}X_{l_1}=X_t\rightarrow X_k$ in $\Mm$. Therefore $t$ has to be $<k$. Moreover since $r-h+1\leq 0$, $l_1$ is $\leq s$ by step 1 of proposition~\ref{Wreduit}. Finally we get $l_1<k$.

Combining substep 1 and substep 2, we can prove that $lj_1j_2\ldots j_s\varphi(i)$ can not be a subword of $w_i$.  Indeed, assume $lj_1j_2\ldots j_s\varphi(i)$ is a subword of $w_i$. There exist $1\leq i_0<i_1<\ldots<i_s<i_{s+1}\leq i$ such that $\varphi(i_0)\varphi(i_1)\ldots\varphi(i_{s+1})=lj_1j_2\ldots j_s\varphi(i)$. In particular, the word $j_1j_2\ldots j_s\varphi(i)$ is a subword of $w_i$ and $1\leq i_1<\ldots<i_s<i_{s+1}\leq i$ is in the set of substep 1. Thus by substep 1, $i_1$ has to be $\leq l_1$. By substep 2, $i_0$ can not exist.

\end{proof}

\subsubsection*{Endomorphism algebra of the cluster-tilting object}

\begin{lema}
 Let $X_i$ and $X_j$ be indecomposables of $\Mm$. We have an isomorphism of vector spaces
$$\Hom_\Lambda(e_{\varphi(j)}\Lambda/\Ii_{w_j},e_{\varphi(i)}\Lambda/\Ii_{w_i})\simeq \bigoplus_{p\geq 0}\Mm(\tau_B^pX_j,X_i).$$
\end{lema}

\begin{proof}
 \textit{Case 1: $j\geq i$}\smallskip\\
By \cite{Bua2} (lemma II.1.14) we have an isomorphism
$$\Hom_\Lambda(e_{\varphi(j)}\Lambda/\Ii_{w_j},e_{\varphi(i)}\Lambda/\Ii_{w_i})\simeq e_{\varphi(i)}\Lambda/\Ii_{w_i}e_{\varphi(j)}.$$
By proposition \ref{imagedeF}, this is isomorphic to the space
$$\bigoplus_{p\geq 0}\Mm(\tau_B^pH_{\varphi(j)},X_i).$$
By definition of the function $\varphi$, there exists some $q\geq 1$ such that $X_j=\tau_B^qH_{\varphi(j)}.$  Thus we can write the sum
$$\bigoplus_{p\geq 0}\Mm(\tau_B^pH_{\varphi(j)},X_i)=\bigoplus_{p=1}^{q}\Mm(\tau_B^{-p}X_j,X_i)\oplus \bigoplus_{p\geq 0}\Mm(\tau_B^pX_j,X_i)$$
Since $j\geq i$, there is no morphism from $\tau_B^{-p}X_j$ to $X_i$ for $p\geq 1$, and the first summand is zero. Therefore we get the result.\medskip\\

\textit{Case 2: $j<i$}\smallskip\\
By \cite{Bua2} (lemma II.1.14) we have an isomorphism
$$\Hom_\Lambda(e_{\varphi(j)}\Lambda/\Ii_{w_j},e_{\varphi(i)}\Lambda/\Ii_{w_i})\simeq e_{\varphi(i)}(\Ii_{\varphi(i)}\ldots\Ii_{\varphi(j+1)}/\Ii_{w_i})e_{\varphi(j)}.$$ By proposition \ref{imagedeF}, this space is a subspace of the space $$\bigoplus_{p\geq 0}\Mm(\tau_B^pH_{\varphi(j)},X_i)\simeq \bigoplus_{p\geq 1}\Mm(\tau_B^{-p}X_j,X_i)\oplus\bigoplus_{p\geq 0}\Mm(\tau_B^pX_j,X_i).$$\medskip\\

\textit{Step 1: If $f$ is a non zero morphism  in $\Mm(\tau_B^{-p}X_j,X_i)$ with $p\geq 1$ then $f$ is not in the space $e_{\varphi(i)}\Ii_{\varphi(i)}\ldots\Ii_{\varphi(j+1)}e_{\varphi(j)}$.}\smallskip\\
Let $X_{l_0}$ be the indecomposable $\tau_B^{-p}X_j$. Since $p\geq 1$ then $l_0$ is $\leq j+1$. The morphism is a sum of composition of the form 
$$\xymatrix{X_{l_0}\ar[r]&X_{l_1}\ar[r] &\cdots\ar[r] & X_{l_r}\ar[r] & X_i}$$ with the $X_{l_k}$ indecomposables. Since $f$ is not zero, we have $j+1\leq l_0 <l_1<\ldots <l_r <i$.  Thus the word $\varphi(l_0)\varphi(l_1)\ldots\varphi(l_r)\varphi(i)$ is a subword of $\varphi(j+1)\varphi(j+2)\ldots\varphi(i)$. Since it holds for each factorization of $f$, the morphism $f$ is not in the space $e_{\varphi(i)}\Ii_{\varphi(i)}\ldots\Ii_{\varphi(j+1)}e_{\varphi(j)}$.\medskip\\

\textit{Step 2: If $f$ is a morphism in $\Mm(\tau_B^{p}X_j,X_i)$ with $p\geq 0$ then $f$ is in the space $e_{\varphi(i)}\Ii_{\varphi(i)}\ldots\Ii_{\varphi(j+1)}e_{\varphi(j)}$.}\smallskip\\
Let $X_{l_0}$ be the indecomposable $\tau^p_BX_j$. Since $p$ is $\geq 0$, we have $l_0\leq j$. Let us show that if $f$ is a composition of irreducible morphisms 
$$\xymatrix{X_{l_0}\ar[r]&X_{l_1}\ar[r] &\cdots\ar[r] & X_{l_r}\ar[r] & X_{l_{r+1}}=X_i}$$ then the word $\varphi(l_0)\varphi(l_1)\cdots\varphi(l_r)\varphi(i)$ is not a subword of  $\varphi(j+1)\varphi(j+2)\ldots\varphi(i)$. 

We have $l_0 <l_1<\cdots<l_r<i$. Since $l_0$ is $<j+1$, and $i$ is $\leq j+1$, there exists $1\leq  k\leq r+1$ such that $l_{k-1}<j+1\leq l_k$. Therefore $\varphi(l_k)\ldots\varphi(l_r)\varphi(i)$ is a subword of $\varphi(j+1)\varphi(j+2)\ldots\varphi(i)$, and the sequence $l_k<l_{k+1}<\cdots <l_r<i$ is the maximal element of the set $$\{j+1\leq i_k<\cdots<i_{r+1}\leq i \quad | \quad \varphi(i_k)=\varphi(l_k),\ldots,\varphi(i_r)=\varphi(l_r),\varphi(i_{r+1})=\varphi(i)\}$$
for the lexicographic order (exactly for the same reasons as in substep 1 of proposition~\ref{imagedeF}).
Now we can prove  exactly as in substep 2 of proposition~\ref{imagedeF} that $\varphi(l_{k-1})$ does not belong to the set $\{\varphi(j+1),\ldots,\varphi(l_k-1)\}$. Thus $\varphi(l_{k-1})\varphi(l_k)\ldots\varphi(l_r)\varphi(i)$ can not be a subword of $\varphi(j+1)\varphi(j+2)\ldots\varphi(i)$.  

Finally, let $f=f_1+f_2$ be a morphism in $$\bigoplus_{p\geq 0}\Mm(\tau_B^pH_{\varphi(j)},X_i)\simeq \bigoplus_{p\geq 1}\Mm(\tau_B^{-p}X_j,X_i)\oplus\bigoplus_{p\geq 0}\Mm(\tau_B^pX_j,X_i).$$ By step 2, $f_2$ is in the space $e_{\varphi(i)}\Ii_{\varphi(i)}\ldots\Ii_{\varphi(j+1)}e_{\varphi(j)}$. By step 1 the morphism $f$ is in  $e_{\varphi(i)}\Ii_{\varphi(i)}\ldots\Ii_{\varphi(j+1)}e_{\varphi(j)}$ if and only if $f-1$ is zero. Thus we get an isomorphism
$$\Hom_\Lambda(e_{\varphi(j)}\Lambda/\Ii_{w_j},e_{\varphi(i)}\Lambda/\Ii_{w_i})\simeq \bigoplus_{p\geq 0}\Mm(\tau_B^pX_j,X_i).$$
\end{proof}

\begin{cora}\label{addTaddA2}
 If $X_i$ and $X_j$ are indecomposables of $\Mmm$, then we have 
$$\Hom_{\underline{Sub}\Lambda/\Ii_w}(e_{\varphi(j)}\Lambda/\Ii_{w_j},e_{\varphi(i)}\Lambda/\Ii_{w_i})\simeq e_{X_j}\tilde{A}e_{X_i}.$$
\end{cora}

\begin{proof}
 The proof is exactly the same as the proof of corollary \ref{addTaddA}.
\end{proof}

\subsubsection*{Triangle equivalence}

\begin{thma}\label{casBIRS}
  The functor $F\circ i_*:\mod \Mmm \rightarrow f.l.\Lambda$ induces an algebraic triangle equivalence between $\Cc_{\Mmm}$ and $\underline{\Sub} \Lambda/\Ii_w$.
\end{thma}

\begin{proof}
By proposition \ref{imagedeF}, the functor $F$ sends the projectives of $\mod \Mm$ onto the summands of the cluster-tilting object $T_w$ of the category $\Sub \Lambda/\Ii_w$. For $i=1,\ldots, n$, the projective $H_i^\wedge$ is sent to the projective-injective $\Lambda/\Ii_w e_i$. Furthermore, by corollary \ref{addTaddA2}, $F\circ i_*$ induces an equivalence between the subcategories $add(A)$ and $add(T_w)$. Thus we can conclude as in the proof of theorem \ref{casGLS}.
 \end{proof}

\subsection{Example} We refer to \cite{Ami2} for more examples.
Let $Q$ be the following quiver: $\xymatrix{1\ar[r]& 2 &3 \ar@<.5ex>[l]\ar@<-.5ex>[l] }$.
The preinjective component of $\mod kQ$ looks as follows:
\smallskip\\

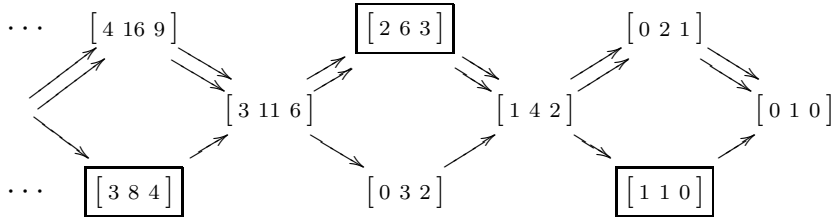
\xymatrix@-1.3pc{\cdots&\left[{\bsm4&16&9\esm}\right]\ar@<.5ex>[dr] \ar@<-.5ex>[dr]&&\fbox{$\left[{\bsm2&6&3\esm}\right]$}\ar@<.5ex>[dr] \ar@<-.5ex>[dr]&&\left[{\bsm0&2&1\esm}\right]\ar@<.5ex>[dr]\ar@<-.5ex>[dr] &\\
\ar@<-.5ex>[ur]\ar@<.5ex>[ur]\ar[dr]&&\left[{\bsm3&11&6\esm}\right]\ar@<.5ex>[ur]\ar@<-.5ex>[ur]\ar[dr]&&\left[{\bsm1&4&2\esm}\right]\ar@<-.5ex>[ur]\ar@<.5ex>[ur]\ar[dr]&&\left[{\bsm0&1&0\esm}\right]\\
\cdots&\fbox{$\left[{\bsm3&8&4\esm}\right]$}\ar[ur]&&\left[{\bsm0&3&2\esm}\right]\ar[ur]&&\fbox{$\left[{\bsm1&1&0\esm}\right]$}\ar[ur]&}

 Here we denote the $kQ$-modules by their dimension vectors in order to lighten the writing. For example the module $\left[{\bsm1&4&2\esm}\right]$ has the following composition series: ${\bsm 2&&2&&2&&2\\&3&&1&&3&\esm}$.

If we mutate the tilting object $\left[{\bsm2&6&3\esm}\right] \oplus\left[{\bsm1&4&2\esm}\right]\oplus\left[{\bsm1&1&0\esm}\right]$ in the direction $\left[{\bsm1&4&2\esm}\right]$, we stay in the preinjective component. We get the tilting object:
$$T=\left[{\bsm2&6&3\esm}\right]\oplus\left[{\bsm3&8&4\esm}\right]\oplus\left[{\bsm1&1&0\esm}\right].$$
The algebra $B=\End_{kQ}(T)$ is a concealed algebra and  is given by the quiver:

$$
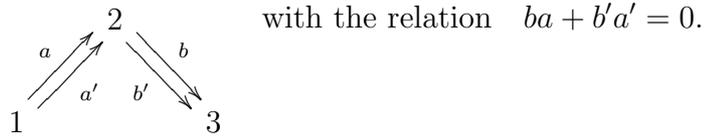
\xymatrix{&2\ar@<.5ex>[dr]^b\ar@<-.5ex>[dr]_{b'}&\\1\ar@<.5ex>[ur]^a\ar@<-.5ex>[ur]_{a'}&&3} \quad \textrm{with the relation}\quad ba+b'a'=0.$$

The functor $R\Hom_{kQ}(T,?)$ yields an equivalence between $\Dd^b(kQ)$
and $\Dd^bB$. Denote by $H$ the image of $D(kQ)$ through this equivalence. This is a
postprojective slice of $\mod B$. Moreover, this equivelence restricts to an equivalence between the category $\Mm =\{X\in \mod B \ | \ \Ext^1_{B}(X,H)=0\}$ and the category $\Mm'=\{X\in \mod kQ \ | \ \Ext^1_{kQ}(T,X)=0\}$. The indecomposable objects of $\Mm'$ are 
$$\left[{\bsm3&8&4\esm}\right], \left[{\bsm2&6&3\esm}\right], \left[{\bsm1&4&2\esm}\right], \left[{\bsm1&1&0\esm}\right], \left[{\bsm0&2&1\esm}\right], \textrm{and} \left[{\bsm0&1&0\esm}\right].$$
The quiver of $\Mm'$ with an admissible ordering is the following:
$$
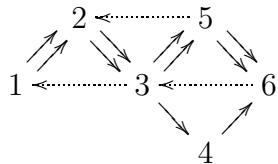
\xymatrix@-1.1pc{&2\ar@<.5ex>[dr] \ar@<-.5ex>[dr]&&5\ar@<.5ex>[dr]\ar@<-.5ex>[dr]\ar@{.>}[ll] &\\
1\ar@<.5ex>[ur]\ar@<-.5ex>[ur]&&3\ar@<-.5ex>[ur]\ar@<.5ex>[ur]\ar[dr]\ar@{.>}[ll]&&6\ar@{.>}[ll].\\
&&&4\ar[ur]&}$$
The dotted arrows represent the Auslander translation $\tau_B$.
The projective indecomposables of $\mod \Mm$ have the following dimension vectors:
$$\left[{\bsm&0&&0&\\1&&0&&0\\&&&0\esm}\right],\quad\left[{\bsm&1&&0&\\2&&0&&0\\&&&0&\esm}\right],\quad\left[{\bsm&2&&0&\\3&&1&&0\\&&&0&\esm}\right],\quad\left[{\bsm&2&&0&\\3&&1&&0\\&&&1&\esm}\right],\quad\left[{\bsm&3&&1&\\4&&2&&0\\&&&0&\esm}\right],\quad\left[{\bsm&6&&2&\\8&&4&&1\\&&&1&\esm}\right]$$

Now let $\Lambda$ be the preprojective associated to the quiver $Q$. The functor $F:\mod \Mm\rightarrow \mod \Lambda$ sends the simples $\Mm$-modules $S_1=\left[{\bsm&0&&0&\\1&&0&&0\\&&&0\esm}\right]$, $S_3=\left[{\bsm&0&&0&\\0&&1&&0\\&&&0\esm}\right]$ and $S_6=\left[{\bsm&0&&0&\\0&&0&&1\\&&&0\esm}\right]$ on the simple $\Lambda$-module $S_2=\left[{\bsm0\\1\\0\esm}\right]$. 
$$\xymatrix@-1.1pc{&2\ar@<.5ex>[dr] \ar@<-.5ex>[dr]\ar@{.}[rr]&&5\ar@{.}[rr]\ar@<.5ex>[dr]\ar@<-.5ex>[dr] &&3\ar@<-.5ex>@{-}[d]\ar@<.5ex>@{-}[d]\\
1\ar@<.5ex>[ur]\ar@<-.5ex>[ur]\ar@{.}[rr]&&3\ar@<-.5ex>[ur]\ar@<.5ex>[ur]\ar[dr]\ar@{.}[rr]&&6\ar@{.}[r]&2\\
&&&4\ar[ur]\ar@{.}[rr]&&1\ar@{-}[u]}.$$
It sends the simple $\Mm$-modules $S_2= \left[{\bsm&1&&0&\\0&&0&&0\\&&&0\esm}\right]$ and $S_5=\left[{\bsm&0&&1&\\0&&0&&0\\&&&0\esm}\right]$ on the simple $\Lambda$-module $S_1=\left[{\bsm1\\0\\0\esm}\right]$, and the simple $\Mm$-module $S_4=\left[{\bsm&0&&0&\\0&&0&&0\\&&&1\esm}\right]$ on the simple $\Lambda$-module $S_3=\left[{\bsm0\\0\\1\esm}\right]$. Since it is exact, it preserves the composition series and then it is easy to compute the image of the indecomposable projective $\Mm$-modules. We get 
$$\left[{\bsm0\\1\\0\esm}\right], \left[{\bsm1\\2\\0\esm}\right], \left[{\bsm2\\4\\0\esm}\right], \left[{\bsm2\\4\\1\esm}\right], \left[{\bsm4\\6\\0\esm}\right] \textrm{and}\left[{\bsm8\\13\\1\esm}\right].$$
 
The projectives of the preprojective algebra associated to $Q$ have the following composition series:
$${\bsm &&&&1&&&&\\&&&&2&&&&\\&&&3&&3&&&\\&&2&&2&&2&&\\3&1&3&3&1&3&3&1&3\\&&\vdots&&\vdots&&\vdots&&\esm},\quad {\bsm&&&&2&&&&\\&&3&&1&&3&&\\&2&&2&&2&&2&\\3&1&3&3&1&3&3&1&3\\&&\vdots&&\vdots&&\vdots&&\esm}, \quad \textrm{and} \quad {\bsm&&&&&3&&&&&\\&&&2&&&&2&&&\\&3&&1&&3&&1&&3&\\2&&2&&2&&2&&2&&2\\&&&&\vdots&&\vdots&&&&\esm}$$

The word $w$ associated with the ordering is $w=232132$. Thus the maximal rigid object of the category $\Sub \Lambda/\Ii_w$ is 
$$R={\bsm2\esm}\oplus{\bsm&3&\\2&&2\esm}\oplus{\bsm&&2&&\\&3&&3&\\2&&2&&2\esm}\oplus{\bsm&&1&&\\&&2&&\\&3&&3&\\2&&2&&2\esm}\oplus{\bsm&&&3&&&\\&&2&&2&&\\&3&&3&&3&\\2&&2&&2&&2\esm}\oplus{\bsm&&&&2&&&&\\&&&3&1&3&&&\\&&2&2&&2&2&&\\&3&3&3&&3&3&3&\\2&2&2&2&&2&2&2&2\esm}.$$ It is easy to check that $R$ is the image by $F$ of the projective indecomposable $\Mm$-modules. The last three summands are the projective-injectives of the Frobenius category $\Sub\Lambda/\Ii_w$.
 This confirms proposition \ref{imagedeF}.

Now take the module $X=1$ in $\Mm$. It corresponds to the module $\left[{\bsm3&8&4\esm}\right]$ in $\mod kQ$. We have the injective resolution in $\mod kQ$:
$$\xymatrix{0\ar[r] & \left[{\bsm3&8&4\esm}\right]\ar[r] &\left[{\bsm0&2&1\esm}\right]^4\oplus\left[{\bsm1&1&0\esm}\right]^3\ar[r]& \left[{\bsm0&1&0\esm}\right]^3\ar[r] & 0}$$
Thus the short exact sequence in $\Mm$: $\xymatrix{0\ar[r] & X\ar[r]& H_0\ar[r] &H_1\ar[r] & 0}$ is the following:
$$\xymatrix{0\ar[r] & 1\ar[r]& 4^3\oplus5^4\ar[r] &6^3\ar[r] & 0}$$
Therefore, the sequence $\xymatrix@-1.1pc{0\ar[r] & X^\wedge\ar[r] & H_0^\wedge\ar[r] & H_1^\wedge\ar[r] & (\tau^{-1}X)^\vee/\proj B\ar[r] & 0}$ in $\mod \Mm$ becomes:
$$\xymatrix@-1.1pc{0\ar[r] & \left[{\bsm&0&&0&\\1&&0&&0\\&&&0&\esm}\right]\ar[r]&\left[{\bsm&2&&0&\\3&&1&&0\\&&&1&\esm}\right]^3\oplus\left[{\bsm&3&&1&\\4&&2&&0\\&&&0&\esm}\right]^4\ar[r]&\left[{\bsm&6&&2&\\8&&4&&1\\&&&1&\esm}\right]^3\ar[r] & \left[{\bsm&0&&2&\\0&&1&&3\\&&&0&\esm}\right]\ar[r] & 0}$$
where $\left[{\bsm&0&&2&\\0&&1&&3\\&&&0&\esm}\right]$ is  the quotient of $(\tau_B^{-1}1)^\vee=3^\vee=\left[{\bsm&0&&2&\\0&&1&&4\\&&&1&\esm}\right]$ by the projectives.
Applying the projection functor we get the exact sequence in $\mod \Lambda$:
$$\xymatrix{0\ar[r] & \left[{\bsm0\\1\\0\esm}\right]\ar[r] &\left[{\bsm2\\4\\1\esm}\right]^3\oplus\left[{\bsm4\\6\\0\esm}\right]^4\ar[r] & \left[{\bsm8\\13\\1\esm}\right]^3\ar[r] &\left[{\bsm2\\4\\0\esm}\right]\ar[r] & 0}$$

The algebra $A$ is the endomorphism algebra of the direct sum of the indecomposables of $\Mmm=\Mm/add H\simeq \Mm'/add D(kQ)$. Thus the 
algebra $A$ is given by the quiver 
$$
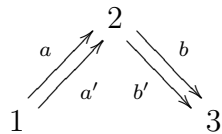
\xymatrix{&2\ar@<.5ex>[dr]^b\ar@<-.5ex>[dr]_{b'}&\\1\ar@<.5ex>[ur]^a\ar@<-.5ex>[ur]_{a'}&&3 } \quad \textrm{and the relation}\quad ba+b'a'=0.$$
By Theorem \ref{homfini2} the cluster category $\Cc_A$ associated with the algebra $A$ is $2$-Calabi-Yau, $\Hom$-finite and $A\in\Cc_A$ is a cluster-tilting object. Moreover by proposition \ref{carquoisdetildeA}, the quiver of the cluster-tilted algebra $\tilde{A}=\End_{\Cc_A}(A)$ has the form:
$$
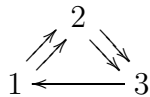
\xymatrix@-1.1pc{&2\ar@<.5ex>[dr]\ar@<-.5ex>[dr]&\\1\ar@<.5ex>[ur]\ar@<-.5ex>[ur]&&3\ar[ll] .}$$

The injective $A$-module $I_1=1^\vee_{|_{\Mmm}}$ has dimension vector $\left[{\bsm&2&\\1&&3\esm}\right]={\bsm3&&3&&3\\&2&&2&\\&&1&&\esm}$. Its image by $i^*$ is the $\Mm$-module  $\left[{\bsm&2&&0&\\1&&3&&0\\&&&0&\esm}\right]$. Its image through $F$ is the same as the image of the $\Mm$-module $\left[{\bsm&0&&2&\\0&&1&&3\\&&&0&\esm}\right]$, indeed we have $F\circ i_*(1^\vee_{|_{\Mmm}})=\left[{\bsm2\\4\\0\esm}\right]$. By the exact sequence above, there is an  isomorphism in $\underline{\Sub} \Lambda/\Ii_w$ between $F\circ i_*( I_1)$ and $F\circ i_*(P_1)[2]$ where $P_1$ is the projective $A$-module with vector dimension $\left[{\bsm&0&\\1&&0\esm}\right]$.

\bibliographystyle{amsalpha}
\bibliography{biblio.bib}

\end{document}